\documentclass[12 pt]{amsart}
\usepackage{graphicx} 
\usepackage{tikz}
\usepackage[margin=0.8in]{geometry}
\usepackage[
backend=biber,
style=alphabetic,
sorting=ynt
]{biblatex}
\usepackage{amsthm}
\usepackage{amsmath}
\usepackage{amssymb}
\usepackage{bbm}
\usepackage{xfrac}
\usepackage{MnSymbol}
\usepackage[english]{babel}
\usepackage{enumitem}
\usepackage{hyperref}
\usepackage{relsize}
\usepackage{quiver}
\usepackage{mathrsfs}
\usepackage{dsfont}
\usepackage{centernot}
\usepackage{mathtools}
\usepackage{stmaryrd}
\usepackage{faktor}
\usepackage{graphicx}
\usepackage{circuitikz}
\usepackage[nice]{nicefrac}
\usepackage{tikz-cd}
\theoremstyle{plain}
\newtheorem{Thm}{Theorem}[section]
\newtheorem*{Thm*}{Theorem}
\newtheorem{Prop}[Thm]{Proposition}
\newtheorem*{Prop*}{Proposition}
\newtheorem{Lemma}[Thm]{Lemma}
\newtheorem*{Lemma*}{Lemma}
\newtheorem{Corollary}[Thm]{Corollary}
\theoremstyle{definition}
\newtheorem{Def}[Thm]{Definition}
\newtheorem*{Not}{Notation}
\newtheorem{notation}[Thm]{Notation}
\newtheorem{Setting}[Thm]{Setting}
\theoremstyle{remark}
\newtheorem*{remark}{Remark}
\DeclareMathOperator{\Frob}{Frob}
\DeclareMathOperator{\Gal}{Gal}
\DeclareMathOperator{\ord}{ord}
\DeclareMathOperator{\lcm}{lcm}
\DeclareMathOperator{\Spec}{Spec}

\DeclareMathOperator{\ScDeg}{ScDeg}
\DeclareMathOperator{\td}{td}
\DeclareMathOperator{\Frac}{Frac}

\DeclareMathOperator{\Aut}{Aut}
\DeclareMathOperator{\id}{id}
\DeclareMathOperator{\Fr}{Fr}

\DeclareMathOperator{\Disc}{Disc}
\DeclareMathOperator{\N}{N}
\DeclareMathOperator{\Cl}{Cl}
\DeclareMathOperator{\codim}{codim}
\DeclareMathOperator{\Sym}{Sym}
\DeclareMathOperator{\lex}{lex}

\DeclareMathOperator{\Sl}{sl}
\DeclareMathOperator{\alg}{alg}
\DeclareMathOperator{\Char}{char}

\DeclareMathOperator{\Norm}{Norm}
\DeclareMathOperator{\Specm}{Specm}
\graphicspath{{}}
\usetikzlibrary{patterns}
\usepackage{pdfpages}

\title{Explicit Error Term in the Chebotarev Density Theorem for Function Fields}
\author{Ido Benbenishty Shpilka}
\date{}

\begin{document}



\maketitle
\begin{abstract}
    We prove an explicit error term in the Chebotarev Density Theorem for finite generically Galois étale morphisms of affine varieties over finite fields, under a mild condition on their size. The error term generalises the known results for curves, and depends only on invariants of the varieties.

    We present several applications to the explicit error term. We apply it to existence problems, and to the Bateman-Horn conjecture over finite fields.
\end{abstract}
\section{Introduction}
The Chebotarev Density Theorem is a fundamental analytic result in algebraic number theory. Over $\mathbb{Q}$, it states that for any Galois extension $K/\mathbb{Q}$ and any automorphism $\sigma \in \Gal(K/\mathbb{Q})$, the set of primes $p \in \mathbb{Z}$ for which there exists 
a prime $P \lhd \mathcal{O}_K$ lying above $p$ such that $\sigma \bmod P \in \Gal\big((\mathcal{O}_K/P)/(\mathbb{Z}/p)\big)$ is defined and acts as the Frobenius automorphism 
has a density of
\[\frac{\lvert\Cl(\sigma)\rvert}{\lvert\Gal(K/\mathbb{Q})\rvert}\, ,\]
where $\Cl(\sigma)$ denotes the conjugacy class of $\sigma$.

For example, if $a,b \in \mathbb{N}$ are coprime, take $K = \mathbb{Q}(\zeta_{a})$, and $\sigma \in \Gal(\mathbb{Q}(\zeta_{a})/\mathbb{Q})$ to be defined by $\sigma(\zeta_{a}) = \zeta_{a}^{b}$. Then, the density of primes associated to $\sigma$ in the notion described above is $\frac{1}{\varphi(a)}$, where $\varphi$ is Euler's totient function. A prime $p$ is associated to $\sigma$ in such a way if and only if $\sigma(\zeta_{a}) = \zeta_{a}^{p}$, which is equivalent to the congruence $p = b \pmod{a}$. Therefore, when applied to cyclotomic extensions, Chebotarev's Density Theorem is equivalent to Dirichlet's Theorem on primes in arithmetic progressions.

Chebotarev's Density Theorem  has further applications to  number theory. 
For example, two Galois extensions $L_1, L_2/\mathbb{Q}$ are equal if and only if the set 
of splitting primes in $L_1/\mathbb{Q}$ coincides with the set of those that split in $L_2/\mathbb{Q}$; 
see \cite[Corollary~VII.13.10]{neukirch}.

Chebotarev's Density Theorem is not limited to $\mathbb{Q}$; it can be formulated over number fields and over function fields over finite fields. In this formulation, we count the prime ideals of the integral closure of $\mathbb{Z}$ or of a polynomial ring, respectively, by their absolute norm.

We now restrict our discussion to Chebotarev's Density Theorem over function fields over a finite field $\mathbb{F}_{q}$. The major step of the proof found in \cite[Section 7.4]{friedjarden} for function fields of dimension 1 is to find an error term for the asymptotics of the number of primes associated to an automorphism $\sigma$ with a residue field of a given cardinality $q^{k}$. For a Galois extension $L/K$ of function fields over $\mathbb{F}_{q}$, the error term in the mentioned proof is given explicitly in terms of $k$, $q$, $\lvert \Cl(\sigma)\rvert$, the genera $g_{K}, g_{L}$ of the fields, and $[L\overline{\mathbb{F}_{q}}:K\overline{\mathbb{F}_{q}}]$, where $\overline{\mathbb{F}_{q}}$ is an algebraic closure of $\mathbb{F}_{q}$, see~\eqref{eq:FJ} below.

We aim to prove an explicit error term for function fields of higher dimensions. Such error terms prove to be important in applications. Before going into the details, we state one such application.

The \emph{factorisation type} of a polynomial $f$ of degree $n$ is the partition $\lambda$ of $n$ induced by the degrees of its irreducible factors.  
We call a polynomial $f\in \mathbb{F}_q[x_1,\ldots, x_d,t]$ \emph{simply ramified in \(t\)} if $\Char \mathbb{F}_{q}>\deg_t f>0$ and for every \(\alpha \in \overline{\mathbb{F}}_q^{\,u}\), the specialisation \(f(\alpha,t)\) is either separable or has a unique multiple root, necessarily of multiplicity \(2\). 
    \begin{Thm}\label{thm:1App}
    Let $f\in \mathbb{F}_{q}[x_{1},\ldots,x_{d},t]$ be a geometrically irreducible polynomial of total degree $N=\deg f>1$ which is simply ramified and monic in $t$. Assume that 
    \[
        q\geq 1336^{2}36^{d}N^{N(2\max\{d,3\}+6)}(\log(N))^{d+3}.
    \]
    Then every partition of $\deg_tf$ occurs as a factorisation type of $f(\alpha,t)$ with $\alpha\in \mathbb{F}_q^{d}$.
    \end{Thm}
    A simpler but weaker lower bound for $q$ for which the theorem holds is
    \[q>CN^{CNd}\, ,\]
    for a sufficiently large $C>0$.
    The simply ramified condition is necessary; for example, take $F=t^N-x$. If $N\mid q-1$, then all irreducible factors of a specialisation have the same degree.

    \subsection{Explicit Chebotarev Theorem for Varieties  over Finite Fields -- Main Results}
 We work in the following setting:
\begin{Setting}\label{notation1}
    1. Let $q$ be a power of a prime. Denote $\mathbb{A}^{u}_{\mathbb{F}_q}=\Spec\mathbb{F}_{q}[x_{1},\ldots,x_{u}]$.
\begin{enumerate}[label=\arabic*.]\addtocounter{enumi}{+1}
    \item Let $X,\ Y$ be normal closed $\mathbb{F}_q$-varieties of dimension $d$, and let $h\colon Y\to X$ be a finite generically Galois morphism. Assume that $X$ is geometrically integral over $\mathbb{F}_{q}$.
    \item Denote $K=\kappa(X)\subseteq \kappa(Y)=L$, $n=[L\cap \overline{\mathbb{F}_{q}}:\mathbb{F}_{q}]$, and $m=[L\overline{\mathbb{F}_{q}}:K\overline{\mathbb{F}_{q}}]$.
    \item Assume $X\subseteq \mathbb{A}^{u}_{\mathbb{F}_{q}},\ Y_{\overline{\mathbb{F}_{q}}}\subseteq \mathbb{A}^{v}_{\overline{\mathbb{F}_{q}}}$ are closed subvarieties of affine spaces. Let $\delta_{X},\ \delta_{Y_{\overline{\mathbb{F}_{q}}}}$ be the degree of $X, Y_{\overline{\mathbb{F}_{q}}}$, respectively (Definition~\ref{Def:2.5}).
    \item Let $l$ denote the number of polynomials $Y_{\overline{\mathbb{F}_{q}}}$ is defined by the vanishing of $l$, and let $D$ denote the maximum of the degrees of these polynomials.
    \item Fix a conjugacy class $\mathcal{C}\subseteq G:=\Gal(L/K)$ and let $r\in \mathbb{N}$ be such that for all $\sigma\in \mathcal{C}$, \[\sigma\vert_{\mathbb{F}_{q^{n}}}=\Frob_{q}^{r}\, .\]
    Such an $r$ exists since $\Gal(\mathbb{F}_{q^n}/\mathbb{F}_{q})$ is abelian, and therefore all elements in $\mathcal{C}$ have the same restriction to $\mathbb{F}_{q^{n}}$.
    \item Let $|X|$ be the set of closed points of $X$ and for every $k\in \mathbb{N}$ let $|X|_k$ the subset of $|X|$ of points of degree $k$. Then, we let 
    \[
        \pi_{\mathcal{C}}(k) = \lvert\{ x\in |X|_k :  x\text{ is unramified in } Y \text{ and } \Fr(x)=\mathcal{C}\}\rvert.
    \]
    \item Let 
    \[
        R_k = \lvert\{ x\in |X|_k: x \text{ ramified in $Y$}\}\rvert.
    \]
\end{enumerate}
\end{Setting}
We want to find an asymptotic formula for $\pi_{\mathcal{C}}(k)$ with an explicit error term. For curves, that is when $d=1$, we have a fully  explicit formula. If $k=r\mod n$, then 
\begin{equation}\label{eq:FJ}
    \left|\dfrac{\pi_{\mathcal{C}}(k)}{q^{k}/k}-\dfrac{\lvert \mathcal{C}\rvert}{m} \right| \leq C_1 q^{-k/2} + C_2q^{-k}
\end{equation}
where $C_1 \leq 2g_Y+4g_Xm+2m$, $C_2 \leq (k+m+1)D$,
and $D = \sum_{k} kR_k$, see  \cite[Theorem 1]{Curves} and \cite[Proposition 7.4.8]{friedjarden}.  Otherwise, if $k\neq r\mod n$, then $\pi_{\mathcal{C}}(k)=0$.

The literature contains several generalisations of \eqref{eq:FJ} to $d>1$, going back at least to Katz, saying that if $Y$ is a closed variety of an affine space, of degree $\delta_{Y}$, $h$ is étale, and that if $k=r\mod n$, then 
\begin{equation}
    \label{eq:CT_Hdime}
    \bigg|\dfrac{\pi_{\mathcal{C}}(k)}{q^{kd}/k}-\frac{\lvert \mathcal{C}\rvert}{m}\bigg|\leq O(q^{-\frac{k}{2}})\, ,
\end{equation}
where the implied constant depends only on $\delta_X,\delta_Y,\deg h,u,v$. See \cite[Theorem~4.1]{Chav} and \cite[Theorem~3]{Entin2} for more precise formulations. 
While these results are uniform, they lack the explicit dependence on those parameters, which may be useful in applications.  

Our main result  is the following explicit version of \eqref{eq:CT_Hdime}. 
\begin{Thm}\label{Thm:1.2}
    Under Setting~\ref{notation1},
    assume that $q>2(d+1)\delta_{X}^{2}$.
    If $k=r\mod n$, then
    \[
        \left|\dfrac{\pi_{\mathcal{C}}(k)}{q^{kd}/k}-\frac{\lvert \mathcal{C}\rvert}{m}\right|  
        \leq  C_1 q^{-\frac{k}{2}} + (C_2+C_3) q^{-k} + C_4q^{-\frac{k}{2}d} + C_5q^{-\frac{k}{2}(d+\frac{1}{2})} + C_6 q^{-\frac{k}{2}(d+1)},
    \]
    where 
    \begin{align*}
    C_1 &= \frac{|\mathcal{C}|}{m}(\delta_{Y_{\overline{\mathbb{F}_{q}}}}-1)(\delta_{Y_{\overline{\mathbb{F}_{q}}}}-2), &
    C_2 &= \frac{6|\mathcal{C}|}{m}\cdot 2^{l}(lD+3)^{v+1},&
    C_3 &= |\mathcal{C}|\,R_k\,kq^{-k(d-1)}, \\[4pt]
    C_4 &= 2|\mathcal{C}|,&
    C_5 &= 2|\mathcal{C}|(\delta_X-1)(\delta_X-2),&
    C_6 &= 10|\mathcal{C}|\delta_X^{\frac{13}{3}}.
\end{align*}
    If $k\neq r\mod n$, then $\pi_{\mathcal{C}}(k)=0$.
\end{Thm}
In the case of $h$ being étale, this generalises \ref{eq:CT_Hdime}.

For $k=1$, we obtain a slightly stronger bound.
\begin{Thm}\label{Thm:1.3}
    Under Setting~\ref{notation1},
    assume that $r=1\mod n$. Then,
    \begin{align*}
    \bigg|\dfrac{\pi_{\mathcal{C}}(1)}{q^{d}}-\frac{\lvert \mathcal{C}\rvert}{m}\bigg|\leq\dfrac{\lvert\mathcal{C}\rvert}{m}\left[(\delta_{Y_{\overline{\mathbb{F}_{q}}}}-1)(\delta_{Y_{\overline{\mathbb{F}_{q}}}}-2)q^{-\frac{1}{2}}+6\cdot 2^{l}(lD+3)^{v+1}q^{-1}\right]+\lvert\mathcal{C}\rvert R_{1}q^{-d}\, .
    \end{align*}
\end{Thm}

We prove Theorems~\ref{Thm:1.2} and \ref{Thm:1.3} in Section~\ref{Section:3}.

As mentioned above, our motivation is derived from a number of applications, which are discussed in Section~\ref{Section:5}. One such application was stated above in Theorem~\ref{thm:1App}. Another line of  applications comes from number theory in function fields, which we now discuss.
\subsection{The Bateman-Horn Conjecture over Finite Fields}
The Bateman-Horn conjecture is a central problem in analytic number theory, concerned with the distribution of prime values of polynomials. Let $\mathcal{F}_{1},\ldots,\mathcal{F}_{l}\in \mathbb{Z}[x]$ be non-associate irreducible polynomials, such that for every prime number $p\in \mathbb{N}$ there exists an $n\in \mathbb{N}$ such that $p\nmid \mathcal{F}_{1}(n),\ldots,\mathcal{F}_{l}(n)$ (for example, the polynomial $x^{2}+x+2$ would be excluded since all of its specialisations are even). The conjecture asserts that there exists a constant $C>0$ dependent only on $\mathcal{F}_{1},\ldots,\mathcal{F}_{l}$ such that
\[\dfrac{\pi_{\mathcal{F}_{1},\ldots,\mathcal{F}_{l}}(N)}{N}=\dfrac{\lvert \{n\in [N]\mid \mathcal{F}_{1}(n),\ldots,\mathcal{F}_{l}(n) \text{ are prime}\}\rvert}{N} \sim C\dfrac{1}{\log(N)^{l}\prod_{i\in [l]} \deg(\mathcal{F}_{i})}\, ,\]
where $[N]=\{n\in \mathbb{N}\mid 1\leq n\leq N\}$.

Heuristically, for a random integer $n\in [N]$, $\mathcal{F}_{i}(n)$ behave like random independent variables of size $\leq a_{i}N^{\deg(\mathcal{F}_{i})}$, where $a_{i}$ are the leading coefficients of $\mathcal{F}_{i}$. Using the prime number theorem, we would get
\begin{align*}
    \dfrac{\pi_{\mathcal{F}_{1},\ldots,\mathcal{F}_{l}}(N)}{N}&=\prod_{i\in [l]} \dfrac{\lvert \{n\in [N]\mid \mathcal{F}_{i}(n)\text{ is prime}\}\rvert}{N}\approx\prod_{i\in [l]} \dfrac{\lvert \{n\in [a_{i}N^{\deg(\mathcal{F}_{i})}]\mid n\text{ is prime}\}\rvert}{a_{i}N^{\deg(\mathcal{F}_{i})}}\\
    &\sim \prod_{i\in [l]}\dfrac{1}{\log(a_{i}N^{\deg(\mathcal{F_{i}})})}\approx \dfrac{1}{\log(N)^{l}\prod_{i\in [l]}\deg(\mathcal{F}_{i})}\, .
\end{align*}
The constant $C$ appears as $\mathcal{F}_{i}(n)$ are coorelated at small primes.
This conjecture generalises many classical conjectures and theorems in number theory. For example, the case $l=1, \mathcal{F}_{1}=ax+b$ where $a,b$ are coprime integers is the prime number theorem for arithmetic progressions, the only proved case of the conjecture. Taking $l=2$ and $\mathcal{F}_{1}=x, \mathcal{F}_{2}=x+2$ would give us a quantitative version of the twin prime conjecture. Lastly, taking $l=1, \mathcal{F}_{1}=x^{2}+1$ would yield Landau's conjecture.

Entin \cite{Entin} following \cite{Pollack,BP,Bary1,SJ,Bary2,BBR} proved an analogue of the Bateman-Horn conjecture over large finite fields. If $\mathcal{F}_{1},\ldots,\mathcal{F}_{l}\in \mathbb{F}_{q}[t][x]$ are non-associate, irreducible, and separable polynomials over $\mathbb{F}_{q}[t]$ that are monic and not constant in $x$, the finite field analogue describes the asymptotic behaviour of
\[
\pi_{q,\mathcal{F}_1,\ldots, \mathcal{F}_l}(d)=\lvert \{f\in \mathbb{F}_{q}[t]\mid \deg f=d, f\text{ is monic},\  \mathcal{F}_{i}(t,f(t)) \text{ is irreducible},\ \forall i\in[l]\}\rvert\, .\]

Entin utilises Chebotarev's Density Theorem in his proof of the analogue result. Plugging in our explicit error term in his proof gives an explicit version of \cite[Theorem 1.1]{Entin}. Recall the definition of the slope of a polynomial from \cite{Entin}:
\begin{Def}
    Let $P=\sum_{j=0}^{d} c_{j}(t)x^{j}\in \mathbb{F}_{q}[t][x]$. The slope of $P$ is defined to be:
    \[\Sl P=\max_{j\in [d]} \frac{\deg c_{d-j}}{j}\]
\end{Def}
Then, the explicit version of \cite[Theorem 1.1]{Entin} goes as follows:
\begin{Thm}\label{BatemanHorn}
    Let $\mathcal{F}_{1},\ldots,\mathcal{F}_{l}\in \mathbb{F}_{q}[t][x]$ be non-associate, irreducible, and separable polynomials over $\mathbb{F}_{q}[t]$, monic $x$, and let
    \[D_{i}=\deg_{x}(\mathcal{F}_{i})>0\, .\]
    Let
    \[d\geq 1+\max_{i\in [l]}\{3,\Sl \mathcal{F}_{i}\}\]
    be a natural number. For all $i\in [l]$, denote $N_{i}=D_{i}d$, $N=\sum_{i\in [l]}N_{i}+D_{i}$, and denote by $\mu_{i}$ the number of irreducible factors into which $\mathcal{F}_{i}(t,x)$ splits over $\overline{\mathbb{F}_{q}}$, so that $\mu_{i}\mid D_{i}$. Set 
    \begin{align*}
    \rho_{i}&=\frac{N_{i}}{\mu_{i}},\ \rho=\prod_{i\in [l]} \rho_{i},\ &
    \gamma&=(\prod_{i\in [l]}(\rho_{i}!)^{\mu_{i}})\lcm(\mu_{1},\ldots,\mu_{l}),\\
    D&=9N^{2}\log(N)\gamma(2\gamma-1)-3, &
    \delta&=DN\gamma\, .
    \end{align*}
    
    Then,
    \begin{align*}
    \lvert \pi_{q,\mathcal{F}_1,\ldots, \mathcal{F}_l}(d)-\rho^{-1}q^{d}\rvert\leq \rho^{-1}q^{d}[(\delta-1)(\delta-2)q^{-\frac{1}{2}}+24(2D+3)^{d+3}q^{-1}]+2(D-1)q^{d-1}\, .
    \end{align*}
    In particular,
    \[\pi_{q,\mathcal F_1,\dots,\mathcal F_l}(d)\sim\rho^{-1} q^{d} \qquad\text{as } q\to\infty.\]
\end{Thm}
\begin{remark}
    In Theorem~\ref{BatemanHorn}, $q$ must grow for the error term to go to zero. If $q>75^{2}\delta^{2d+6}$, then the error term is at most
    \[\dfrac{1}{3}\rho^{-1}q^{d}\, ,\]
    and therefore
\[
\pi_{q,\mathcal{F}_1,\dots,\mathcal{F}_l}(d)
>
\frac{2}{3}\rho^{-1}q^{d}.
\]
In particular,
\[
\pi_{q,\mathcal{F}_1,\dots,\mathcal{F}_l}(d)>0.
\]

Theorem~\ref{BatemanHorn} should be compared with the much deeper results of
Sawin and Shusterman~\cite{SS1,SS2}, where $q$ is fixed and
the degree $d$ tends to infinity. Their results apply to the cases of
two linear polynomials and of one quadratic polynomial, under additional
technical assumptions on $q$ (in particular, $q$ cannot be prime in those
settings). Thus, even in these special cases, Theorem~\ref{BatemanHorn} yields some new
information by providing a completely explicit and uniform bound for
all sufficiently large finite fields.
\end{remark}
\subsection{On The Proof}
To prove Theorem~\ref{Thm:1.2}, we apply the so called field crossing argument to reduce the problem to the case where $L/K$ is a scalar extension and $k=1$. Then, we reduce to estimating the number of rational points on varieties, which is classically done by Lang-Weil estimates. This is an adaptation of the proof presented in \cite[Section 7.4]{friedjarden}.

When $d=1$, it follows from the Riemann Hypothesis for Function Fields that $\lvert X(\mathbb{F}_{q})-q\rvert\leq 2g_{X}\sqrt{q}+1$, leading to ~\ref{eq:FJ}. A theorem of Cafure and Matera (\cite[Theorem 7.1]{cafurematera}) allows us to obtain a similar approximation for varieties of higher dimensions. This theorem serves as the main tool for our proof. 

\section{Background in Arithmetic Geometry}
We use the language of arithmetic geometry in this thesis. In particular, we assume that the reader is accustomed with the basic notions of schemes. In particular, with the definition of a scheme $S$ as a locally ringed space \cite[\href{https://stacks.math.columbia.edu/tag/01IJ}{Tag 01IJ}]{stacks-project}, stalks \cite[\href{https://stacks.math.columbia.edu/tag/0078}{0078}]{stacks-project} morphisms locally of finite type \cite[\href{https://stacks.math.columbia.edu/tag/01T1}{01T1}]{stacks-project}, finite morphisms \cite[\href{https://stacks.math.columbia.edu/tag/01WH}{01WH}]{stacks-project}, dominant morphisms \cite[\href{https://stacks.math.columbia.edu/tag/01RJ}{01RJ}]{stacks-project}, separated schemes \cite[\href{https://stacks.math.columbia.edu/tag/01KK}{01KK}]{stacks-project}, integral schemes \cite[\href{https://stacks.math.columbia.edu/tag/01OK}{01OK}]{stacks-project}, geometrically integral scheme \cite[\href{https://stacks.math.columbia.edu/tag/020H}{020H}]{stacks-project}, irreducible scheme \cite[\href{https://stacks.math.columbia.edu/tag/004V}{004V}]{stacks-project}, absolutely irreducible scheme \cite[\href{https://stacks.math.columbia.edu/tag/0365}{0365}]{stacks-project}, normal scheme \cite[\href{https://stacks.math.columbia.edu/tag/033I}{033I}]{stacks-project}, scalar extension \cite[\href{https://stacks.math.columbia.edu/tag/01JX}{01JX}]{stacks-project}.

We recall definitions and basic properties of varieties over fields, as those  may vary in different sources.
Let $E$ be a field. For every extensions $F/E$, denote 
\[
    F_{\alg}:=\{\alpha\in F\mid \alpha \text{ is algebraic over } E\}\, .
\]
\begin{Def}[Variety Over a Field]
    A variety over $E$, or an $E$-variety, is a separated, integral scheme $X$ with a morphism of finite type $X\to \Spec E$.
\end{Def}
\begin{Def}[Affine Variety Over a Field]
    An affine $E$-variety is an $E$-variety $X$ given with a closed embedding $X\hookrightarrow \mathbb{A}_{E}^{n}$ for some $n\in \mathbb{N}$.
\end{Def}

\begin{Def}[Generically Galois Morphism]
    A dominant morphism $h\colon Y\to X$ of $E$-varieties is said to be generically Galois if and only if the associated field extension $\kappa(Y)/\kappa(X)$ is Galois.

    In particular, $\dim X=\dim Y$.
\end{Def}
For the rest of the section, we work in the following setting.
\begin{Setting}\label{Setting:first}
    Let $X,Y$ be normal affine $E$-varieties, with $X$ geometrically integral. Let
\[\begin{tikzcd}
	{Y} & {X} & {\Spec E}
	\arrow["h", from=1-1, to=1-2]
	\arrow[from=1-2, to=1-3]
\end{tikzcd}\]
be a finite generically Galois morphism.
Denote by
\[K=\kappa(X)\subseteq L=\kappa(Y)\]
the function fields of $X,\ Y$ respectively, where the extension is induced by $h$. Then, $L/K$ is a finite Galois extension. Denote:
\[G=\Gal(L/K)\, .\]
\end{Setting}
\begin{Def}
    Let $x\in X$. Denote by $m_{x}\lhd \mathcal{O}_{X,x}$ the unique maximal ideal of the local ring $\mathcal{O}_{X,x}$.
    We define
    \[\kappa(x):=\mathcal{O}_{X,x}/m_{x}\]
    to be the residue field at $x$.
\end{Def}
\begin{remark}\label{remark:finiteness}
    Let $x\in U\subseteq X$ be an affine neighbourhood of $x$, $A=\mathcal{O}_{X}(U)$, and $p\lhd A$ the prime ideal corresponding to $x$. Then:
    \[\kappa(x)\cong A_{p}/pA_{p}\cong \Frac(A/p)\]
    This isomorphism is induced from the map $A=\mathcal{O}_{X}(U)\to \mathcal{O}_{X,x}$. Hence, from now on we identify these fields as the same.
    
    Assume further that $x$ is closed. Then, $\kappa(x)\cong A/p$. Since $A$ is a finitely generated $E$-algebra, it follows that $\kappa(x)$ is a finite extension of $E$. In particular, if $E$ is a finite field, then $\kappa(x)$ is also a finite field.
\end{remark}
\begin{Def}
    Let $x\in X$. Denote by
    \[\mathcal{O}_{Y,x}\subseteq L\]
    the integral closure of $\mathcal{O}_{X,x}$ in $L$.
    
    Then, $h^{-1}(x)=\Specm \mathcal{O}_{Y,x}$, naturally, where $y\in h^{-1}(x)$ corresponds to a maximal ideal $m_{y}\lhd \mathcal{O}_{Y,x}$ satisfying $(\mathcal{O}_{Y,x})_{m_{y}}=\mathcal{O}_{Y,y}$. 
\end{Def}
\begin{remark}
    We list some useful properties of these ideals:
    \begin{enumerate}[label=\arabic*.]
        \item every maximal ideal of $\mathcal{O}_{Y,x}$ is equal to some $m_{y}$.
        \item For each $y\in h^{-1}(x)$, $m_{y}$ lies above $m_{x}$.
    \end{enumerate}
\end{remark}
A key notion in this work is ramification, which we now discuss.
\begin{Prop}\label{Prop:unram}
    Assume that $E$ is a perfect field, and let $x\in X$ be a point in $X$. Denote the following conditions:
    \begin{enumerate}[label=\arabic*.]
        \item $m_{x}\mathcal{O}_{Y,x}=\bigcap_{y\in h^{-1}(x)} m_{y}$.
        \item For all $y\in h^{-1}(x)$, $m_{x}\mathcal{O}_{Y,y}=m_{y}\mathcal{O}_{Y,y}$.
        \item $\mathcal{O}_{Y,x}/m_{x}\mathcal{O}_{Y,x}$ is isomorphic to a product of separable field extensions of $\kappa(x)$.
    \end{enumerate}
    Then, $3\Rightarrow 1\Leftrightarrow 2$, and in the case that $x$ is closed, the conditions are equivalent.
\end{Prop}
\begin{proof}
    $\underline{1\Rightarrow 2}:$
    Localise $m_{x}\mathcal{O}_{Y,x}=\bigcap_{y\in h^{-1}(x)} m_{y}$ at $y\in h^{-1}(x)$. Since for every $y\neq y'\in h^{-1}(x)$, $m_{y'}\nsubseteq m_{y}$, then $m_{y'}\mathcal{O}_{Y,y}=\mathcal{O}_{Y,y}$ and we get $m_{x}\mathcal{O}_{Y,y}=m_{y}\mathcal{O}_{Y,y}$.

    $\underline{2\Rightarrow 1}:$
    $m_{x}\mathcal{O}_{Y,x}\subseteq \bigcap_{y\in h^{-1}(x)} m_{y}$. These ideals are equal, since their localisations at maximal ideals are equal, and equality is a local property.
    
    $\underline{3\Rightarrow 1}:$
    A product of fields has no nilpotent elements and so $m_{x}\mathcal{O}_{Y,x}$ is radical. Thus, $m_{x}\mathcal{O}_{Y,x}$ is the intersection of all primes in $\mathcal{O}_{Y,x}$ that lie above $m_{x}$. The primes that lie above $m_{x}$ are maximal, so $m_{x}\mathcal{O}_{Y,x}$ is the intersection of all maximals containing it. In other words, $m_{x}\mathcal{O}_{Y,x}=\bigcap_{y\in h^{-1}(x)} m_{y}$.

    Now assume $x$ is closed.

    $\underline{1\Rightarrow 3}:$
    Each $m_{y}$ is maximal, so from the Chinese Remainder Theorem,
    \[\mathcal{O}_{Y,x}/m_{x}\mathcal{O}_{Y,x}\cong \prod_{y\in h^{-1}(x)} \mathcal{O}_{Y,x}/m_{y}=\prod_{y\in h^{-1}(x)} \kappa(y)\]
    The field $\kappa(x)$ is perfect as a finite extension of $E$. Hence, $\kappa(y)/\kappa(x)$ are finite, separable extensions.    
\end{proof}
\begin{Def}[Unramified Point]
    A point $x\in X$ is unramified in $Y$ if and only if it satisfies one of the first two equivalent conditions of Proposition~\ref{Prop:unram}.
\end{Def}
\begin{Def}[Unramified Morphism]
    A morphism $f\colon Y\to X$ is said to be unramified if and only if every $x\in \lvert X\rvert$ is unramified in $Y$.
\end{Def}
\begin{Prop}\label{Cor:transitive}
    Assume $h=g\circ f$, where
    \[\begin{tikzcd}
	{Y} & {Z} & {X}
	\arrow["f", from=1-1, to=1-2]
	\arrow["g", from=1-2, to=1-3]
\end{tikzcd}\]
    are morphisms of $E$-varieties. Let $x\in X$. Then:
    \begin{enumerate}[label=\arabic*.]
        \item If $x$ is unramified in $Y$, then every $z\in g^{-1}(x)$ is unramified in $Y$.
        \item If $x$ is unramified in $Z$, and every $z\in g^{-1}(x)$ is unramified in $Y$, then $x$ is unramified in $Y$.
        \item Assume that there exists a finite extension $F/E$ such that $Y$ is an irreducible component of $Z_{F}$. If $x$ is closed and is unramified in $Y$, it is also unramified in $Z$.
    \end{enumerate} 
\end{Prop}
\begin{proof}
   Assume $x$ is unramified in $Y$, and let $z\in g^{-1}(x)$. Then, for all $y\in f^{-1}(z)$, $m_{x}\mathcal{O}_{Y,y}=m_{y}\mathcal{O}_{Y,y}$. Therefore,
    \[m_{y}\mathcal{O}_{Y,y}=m_{x}\mathcal{O}_{Y,y}\subseteq m_{z}\mathcal{O}_{Y,y}\subseteq m_{y}\mathcal{O}_{Y,y}\, ,\]
    and using the second criterion for unramification we conclude that $z$ is unramified in $Y$.

    For $2$, let $y\in h^{-1}(x)$ and $z=f(y)$. Then,
    \[m_{x}\mathcal{O}_{Y,y}=m_{x}\mathcal{O}_{Z,z}\mathcal{O}_{Y,y}=m_{z}\mathcal{O}_{Y,y}=m_{y}\mathcal{O}_{Y,y}\, .\]
    Thus, from the second criterion for unramification, $x$ is unramified in $Y$. 

    Lastly, assume that $Y$ is an irreducible component of $Z_{F}$ and that $x\in X$ is closed and is unramified in $Y$. Then, $\mathcal{O}_{Y,x}=\mathcal{O}_{Z,x}\otimes_{\kappa(Z)_{\alg}} F$, and therefore
    \[\mathcal{O}_{Y,x}/m_{x}\mathcal{O}_{Y,x}\cong \mathcal{O}_{Z,x}/m_{x}\mathcal{O}_{Z,x}\otimes_{\kappa(Z)_{\alg}} F\, .\]
    Thus, from \cite[Theorem 1.3.1]{Milne}, if $\mathcal{O}_{Y,x}/m_{x}\mathcal{O}_{Y,x}$ is isomorphic to a product of separable field extensions of $\kappa(x)$, so is $\mathcal{O}_{Z,x}/m_{x}\mathcal{O}_{Z,x}$, and $x$ is unramified in $Z$.
\end{proof}
    Separable scalar extensions are unramified:
\begin{Prop}\label{Prop:5}
     Let $F/E$ be a separable extension. Then, every $x\in \lvert X\rvert$ is unramified in $X_{F}\to X$.
\end{Prop}
\begin{proof}
     Let $x\in \lvert X\rvert$, and denote by $m_{x}\lhd \mathcal{O}_{X,x}$ the unique maximal ideal of the local ring $\mathcal{O}_{X,x}$. We need to show that in the extension $\mathcal{O}_{X,x}\subseteq \mathcal{O}_{X,x}F$,
     \[\mathcal{O}_{X,x}F/m_{x}\mathcal{O}_{X,x}F\]
     is isomorphic to a product of separable extensions of $\kappa(x)$.
    
     There exists a surjective homomorphism
     \[\kappa(x)\otimes_{E}F\to \mathcal{O}_{X,x}F/m_{x}\mathcal{O}_{X,x}F\, ,\]
     so it is enough to show that $\kappa(x)\otimes_{E}F$ is isomorphic to a product of separable extensions of $\kappa(x)$. This is immediate: write $F\cong E[y]/(f)$ for some separable $f\in E[y]$ to get:
     \[\kappa(x)\otimes_{E}F\cong \kappa(x)\otimes_{E}E[y]/(f)\cong \kappa(x)[y]/(f)\, .\]
     Decomposing $f$ into prime polynomials over $\kappa(x)$ gives the desired result.
\end{proof}
A key notion in this work is the Frobenius conjugacy class. We take a general point of view and define the decomposition and inertia groups, residue degrees, and ramification indices of points.
\begin{Def}
    Let $x\in \lvert X\rvert$, and let $y\in h^{-1}(x)$. Define:
    \begin{enumerate}[label=\arabic*.]
        \item \[D_{y/x}:=\{\sigma\in G\mid \sigma[m_{y}]=m_{y}\}\]
        to be the decomposition group at $y$.
        \item \[I_{y/x}=\{\sigma\in G\mid \forall{a\in \mathcal{O}_{Y,x}}, \sigma(a)-a\in m_{y}\}\lhd D_{y/x} \]
        to be the inertia group at $y$.
        \item The action of $D_{y/x}$ on $\mathcal{O}_{Y,y}$ induces a map
        \[\Psi_{y/x}\colon D_{y/x}\longrightarrow \Gal(\kappa(y)/\kappa(x))\]
        given by $\Psi_{y/x}(\sigma)(b+m_{y})=\sigma(b)+m_{y}$.

        Notice that $I_{y/x}=\ker \Psi_{y/x}$.
        \item The ramification index and the residue degree at $y$ are defined as
        \[e_{y/x}=\lvert I_{y/x}\rvert,\ f_{y/x}=[\kappa(y):\kappa(x)]\]
        respectively.
        Finally, denote
        \[g_{x}=\lvert h^{-1}(x)\rvert\, .\]
    \end{enumerate}
    
\end{Def}
\begin{Prop}
    Assume that $h=g\circ f$, where
    \[\begin{tikzcd}
	{Y} & {Z} & {X}
	\arrow["{f}", from=1-1, to=1-2]
	\arrow["h", curve={height=-24pt}, from=1-1, to=1-3]
	\arrow["g", from=1-2, to=1-3]
\end{tikzcd}\, ,\]
are morphisms of $E$-varieties. Then, for every closed point $x\in Y$, if $y\in h^{-1}(x)$, and $z=f(y)$, then:
\begin{enumerate}[label=\arabic*.]
    \item $D_{y/z}=D_{y/x}\cap \Gal(L/\kappa(Z))$.
    \item $\Psi_{y/x}|_{D_{y/z}}=\Psi_{y/z}$.
    \item $I_{y/z}=I_{y/x}\cap \Gal(L/\kappa(Z))$.
\end{enumerate}
\end{Prop}
\begin{proof}
    Denote by $m_{y/z}\lhd \mathcal{O}_{Y,z}$ the ideal satisfying $(\mathcal{O}_{Y,z})_{m_{y/z}}=\mathcal{O}_{Y,y}$. Then, $m_{y}\mathcal{O}_{Y,y}=m_{y}(\mathcal{O}_{Y,x})_{m_{y}}$ and $m_{y/z}\mathcal{O}_{Y,y}=m_{y/z}(\mathcal{O}_{Y,z})_{m_{y/z}}$ are both maximal ideals of a local ring. Therefore,
    \[m_{y}(\mathcal{O}_{Y,x})_{m_{y}}=m_{y/z}(\mathcal{O}_{Y,z})_{m_{y/z}}\, .\]
    Thus,
    \begin{align*}
        D_{y/z}&=\{\sigma\in \Gal(L/\kappa(Z))\mid \sigma[m_{y/z}]=m_{y/z}\}\\
        &=\{\sigma\in \Gal(L/\kappa(Z))\mid \sigma[m_{y/z}(\mathcal{O}_{Y,z})_{m_{y/z}}]=m_{y/z}(\mathcal{O}_{Y,z})_{m_{y/z}}\}\\
        &=\{\sigma\in \Gal(L/\kappa(Z))\mid \sigma[m_{y}(\mathcal{O}_{Y,x})_{m_{y}}]=m_{y}(\mathcal{O}_{Y,x})_{m_{y}}\}\\
        &=\{\sigma\in \Gal(L/\kappa(Z))\mid \sigma[m_{y}]=m_{y}\}\\
        &=D_{y/x}\cap \Gal(L/\kappa(Z))\, .
    \end{align*}
    This concludes $1$.
    For $2$, let $\sigma\in D_{y/z}$, and let $b\in \mathcal{O}_{Y,x}$. We identify $\mathcal{O}_{Y,x}/m_{y}\cong \mathcal{O}_{Y,z}/m_{y/z}$ through the inclusion $\mathcal{O}_{Y,x}\subseteq \mathcal{O}_{Y,z}$ to get:
    \[\Psi_{y/x}(\sigma)(b+m_{y})=\sigma(b)+m_{y}=\sigma(b)+m_{y/z}=\Psi_{y/z}(\sigma)(b+m_{y/z})\, ,\]
    as desired.
    Condition $3$ follows from comparing the kernels of $\Psi_{y/x}\vert_{D_{y/z}}$ and $\Psi_{y/z}$.
\end{proof}
\begin{Prop}\label{Prop:3}
    Let $x\in \lvert X\rvert$. Then, 
    \begin{enumerate}[label=\arabic*.]
        \item The group $G$ acts transitively on \[\{m_{y}\mid y\in h^{-1}(x)\}\, .\]
        \item Let $y\in h^{-1}(x)$. Then, $\Psi_{y/x}$ is surjective. In particular, \[\lvert D_{y/x}\rvert=f_{y/x}e_{y/x}\, .\]
    \end{enumerate}
\end{Prop}
\begin{proof}
    The assertions are proved in~\cite[Chapter 1, Proposition 9.1, Proposition 9.4]{neukirch}.
\end{proof}
\begin{Corollary}\label{lemma:7}
    Let $x\in \lvert X\rvert$. Then,
    \begin{enumerate}[label=\arabic*.]
        \item $\forall{y,y'\in h^{-1}(x)}$,
        \[D_{y/x}\sim D_{y'/x},\ I_{y/x}\sim I_{y'/x},\ f_{y/x}=f_{y'/x}\, ,\]
        where $\sim$ denotes conjugation in $G$.
        \item $\forall{y\in h^{-1}(x)},$
        \[[L:K]=e_{y/x}f_{y/x}g_{x}\, .\]
    \end{enumerate}
\end{Corollary}
\begin{proof}
    Since $G$ acts transitively on $\{m_{y}\mid y\in h^{-1}(x)\}$, it follows that for every $y,y'\in h^{-1}(x)$, there exists a $\sigma\in G$ for which $\sigma [m_{y}]=m_{y'}$.
    
    Therefore,
    \[\sigma D_{y/x}\sigma^{-1}=D_{y'/x},\ \sigma I_{y}\sigma^{-1}=I_{y'/x}\, .\]
    Then, by Proposition~\ref{Prop:3}, $f_{y/x}=f_{y'/x}$.
    
    For $2$, let $y\in h^{-1}(x)$. It suffices to show that
    \[[L:K]=\lvert D_{y/x}\rvert g_{x}\, ,\]
    which follows from the orbit-stabiliser theorem.
\end{proof}
\begin{notation}
    By Corollary~\ref{lemma:7}, $f_{y/x}$ and $e_{y/x}$ are independent of the choice of $y$, so we simplify and write $f_{x}$ and $e_{x}$, respectively.
\end{notation}
\begin{Corollary}\label{lemma:3.6}
    Let
    \[g\colon Z\to X\]
    be a finite separable morphism of normal affine $E$-varieties. Denote $F=\kappa(Z)$, and let $x\in \lvert X\rvert$.
    Then,
    \[\lvert g^{-1}(x)\rvert\leq [F:K]\, .\]
\end{Corollary}
\begin{proof}
    Let $M$ be the Galois closure of $F/K$, and let $V$ be the normalisation of $Z$ in $M$. Let
    \[\begin{tikzcd}
	{V} & {Z} & X
	\arrow["{f}", from=1-1, to=1-2]
	\arrow["g", from=1-2, to=1-3]
\end{tikzcd}\]
be the corresponding morphism.

    By Corollary~\ref{lemma:7}, we have
    \begin{align*}
        [M:F][F:K]=[M:K]=\sum_{v\in (gf)^{-1}(x)} e_{v/x}f_{v/x}&=\sum_{z\in g^{-1}(x)} \sum_{v\in f^{-1}(z)} e_{v/x}f_{v/z}f_{z/x}\\
        &\geq \sum_{z\in g^{-1}(x)} f_{z/x}\sum_{v\in f^{-1}(z)} e_{v/z}f_{v/z}=[M:F]\lvert g^{-1}(x)\rvert\, .
    \end{align*}
    Therefore, $\lvert f^{-1}(x)\rvert\leq [K':K]$.
\end{proof}
\begin{Corollary}\label{Cor:e=1}
    If $x\in \lvert X\rvert$ is unramified in $Y$, then $e_{x}=1$.
\end{Corollary}
\begin{proof}
    Let $x\in X$ be unramified and $y\in h^{-1}(x)$. Let $Z$ be the normalisation of $X$ in $L^{I_{y/x}}$, and $f\colon Y\to Z,\ g\colon Z\to X$ the corresponding maps. Denote $z=f(y)$, and $m_{z}=m_{y}\cap \mathcal{O}_{Z,x}$:
    \[\begin{tikzcd}
	L & {\mathcal{O}_{Y,x}} & {m_{y}} \\
	{L^{I_{y/x}}} & \mathcal{O}_{Z,x} & {m_{z}} \\
	K & {\mathcal{O}_{X,x}} & {m_{x}}
	\arrow[no head, from=1-1, to=2-1]
	\arrow[no head, from=1-2, to=2-2]
	\arrow[no head, from=1-3, to=2-3]
	\arrow[no head, from=2-1, to=3-1]
	\arrow[no head, from=2-2, to=3-2]
	\arrow[no head, from=2-3, to=3-3]
    \end{tikzcd}\]
    Corollary~\ref{lemma:7} shows that
    \[e_{z}f_{z}g_{z}=[L:L^{I_{y/x}}]=e_{x}\, ,\]
    Since $I_{y/z}=I_{y/x}\cap \Gal(L/L^{I_{y/x}})=I_{y/x}$, we conclude that 
    \[e_{z}=e_{z}f_{z}g_{z}\, .\]
    Therefore, $g_{z}=f_{z}=1$. By Proposition~\ref{Cor:transitive}, $z$ is unramified in $Y$, and using condition $1$ of Proposition~\ref{Prop:unram}, $m_{z}\mathcal{O}_{Y,x}=m_{y}$. Since $f_{z}=1$, it follows that $\mathcal{O}_{Y,x}/m_{y}=\mathcal{O}_{Z,x}/m_{z}$. Therefore, \[\mathcal{O}_{Y,x}=\mathcal{O}_{Z,x}+m_{y}\mathcal{O}_{Z,x}=\mathcal{O}_{Z,x}+m_{z}\mathcal{O}_{Z,x}\, .\]
    Therefore,
    \[\mathcal{O}_{Y,x}/\mathcal{O}_{Z,x}=m_{z}(\mathcal{O}_{Y,x}/\mathcal{O}_{Z,x})\]
    as $\mathcal{O}_{Z,x}$-modules. Since $h$ is a finite morphism, $\mathcal{O}_{Y,x}$ is a finitely generated $\mathcal{O}_{X,x}$-module, and in particular, it is a finitely generated $\mathcal{O}_{Z,x}$-module. Nakayama's lemma implies that $\mathcal{O}_{Y,x}=\mathcal{O}_{Z,x}$. Thus, $I_{y/x}=1$, as wanted.
\end{proof}
From now on, we assume that $E=\mathbb{F}_{q}$ is a finite field. In particulat, $E$ is perfect.
\begin{Def}[Scalar Degree]
    Let $Z$ be an $\mathbb{F}_{q}$-variety. Define the scalar degree to be:
    \[\ScDeg(Z)=\ScDeg_{\mathbb{F}_{q}}(Z):=[\kappa(Z)_{\alg}:\mathbb{F}_{q}]\, .\]
\end{Def}
\begin{Prop}
    An $\mathbb{F}_{q}$-variety $Z$ is geometrically integral over $\mathbb{F}_{q}$ if and only if $\ScDeg_{K}(Z)=1$.
\end{Prop}
\begin{proof}
    $\underline{\Rightarrow}:$ Assume that $Z$ is geometrically integral over $\mathbb{F}_{q}$. Then, $Z_{\overline{\mathbb{F}_{q}}}$ is irreducible, meaning that $\kappa(Z_{\overline{\mathbb{F}_{q}}})=\kappa(Z)\otimes_{\mathbb{F}_{q}} \overline{\mathbb{F}_{q}}$ is a field. Therefore, $\kappa(Z)$ and $\overline{\mathbb{F}_{q}}$ are linearly independent over $\mathbb{F}_{q}$, and in particular, $\ScDeg(Z)=1$.

    $\underline{\Leftarrow}:$ Assume that $\ScDeg(Z)=1$. Then, $\kappa(Z)\otimes_{\mathbb{F}_{q}} \overline{\mathbb{F}_{q}}$ is a field, meaning that $Z_{\overline{\mathbb{F}_{q}}}$ is irreducible.
    Using ~\cite[\href{https://stacks.math.columbia.edu/tag/020I}{020I}]{stacks-project}, we conclude that $Z$ is geometrically integral.
\end{proof}
\begin{Prop}
    Let $Z$ be a normal $\mathbb{F}_{q}$-variety. Then, $\ScDeg_{\mathbb{F}_{q}}(Z)$ is the maximal $n\in \mathbb{N}$ such that $V\to \Spec \mathbb{F}_{q}$ factors through
    \[\begin{tikzcd}
	{Z} & {\Spec \mathbb{F}_{q^{n}}} & {\Spec \mathbb{F}_{q}}
	\arrow[from=1-1, to=1-2]
	\arrow[from=1-2, to=1-3]\, .
\end{tikzcd}\]
\end{Prop}
\begin{proof}
    Let $n\in \mathbb{N}$ be the maximal such natural number. Then, we can embed $\mathbb{F}_{q^{n}}\subseteq \kappa(Z)$, meaning that
    \[\ScDeg(Z)=[\kappa(Z)_{\alg}:\mathbb{F}_{q}]=[\kappa(Z)_{\alg}:\mathbb{F}_{q^{n}}][\mathbb{F}_{q^{n}}:\mathbb{F}_{q}]=n[\kappa(Z)_{\alg}:\mathbb{F}_{q^{n}}]\geq n\, .\]
    Additionally, $\mathbb{F}_{q^{\ScDeg(Z)}}\subseteq \kappa(Z)$, and so if $U\subseteq Z$ is an open affine set, $\mathbb{F}_{q^{\ScDeg(Z)}}\subseteq \mathcal{O}_{Z}(U)$. We therefore have an induced morphism
    \[Z\to \Spec \mathbb{F}_{q^{\ScDeg(Z)}}\, ,\]
    as wanted.
\end{proof}
\begin{remark}
    Since $\mathbb{F}_{q^{\ScDeg(Z)}}\subseteq \kappa(Z)$, then for all $z\in Z$, $\mathbb{F}_{q^{\ScDeg(Z)}}\subseteq \Frac(\mathcal{O}_{Z,z})$. Thus, $\mathbb{F}_{q^{\ScDeg(Z)}}\subseteq \mathcal{O}_{Z,z}$, and we can embed $\mathbb{F}_{q^{\ScDeg(Z)}}\subseteq \mathcal{O}_{Z,z}/m_{z}=\kappa(z)$.
    
    Additionally, as explained in remark~\ref{remark:finiteness}, for every closed point $z\in Z$, $\kappa(z)$ is a finite field.
\end{remark}
\begin{Def}[Norm and Degree]\label{Def:3.4}
Let $Z$ be an $\mathbb{F}_{q}$-variety, and let $z\in \lvert Z\rvert$.
We define the absolute norm of $z$ to be
    \[\N z:=\lvert \kappa(z)\rvert\, ,\]
and the degree of $z$ to be
    \[\deg(z)=\deg_Z(z):=\log_{q^{\ScDeg(Z)}}(\N z)=[\kappa(z):\mathbb{F}_{q^{\ScDeg(Z)}}]\, .\]
\end{Def}
We add to Setting~\ref{Setting:first}:
\begin{notation}\label{Setting}
    Denote
\[n=\ScDeg(Y),\ m=[L:K\mathbb{F}_{q^{n}}],\ d=\dim(X) .\]
\end{notation}
\begin{Corollary}\label{Cor:2}
   Let $x\in \lvert X\rvert$. The following are equivalent:
   \begin{enumerate}[label=\arabic*.]
       \item $e_{x}=1$.
       \item $[L:K]=f_{x}g_{x}$.
       \item For every ${y\in h^{-1}(x)}$, we have $[L:K]=\dfrac{n\deg_{Y}(y)}{\deg_{X}(x)}g_{x}$.
   \end{enumerate}
\end{Corollary}
\begin{proof}
    $(1\Leftrightarrow 2):$
    Immediate from Corollary~\ref{lemma:7}.
    
    $(2\Leftrightarrow 3):$
    Follows from
    \[f_{x}g_{x}=f_{y/x}g_{x}=\dfrac{[\kappa(y):\mathbb{F}_{q^{n}}][\mathbb{F}_{q^{n}}:\mathbb{F}_{q}]}{[\kappa(x):\mathbb{F}_{q}]}g_{x}=\dfrac{n\deg_{Y}(y)}{\deg_{X}(x)}g_{x}\, .\]
\end{proof}

Let $x\in \lvert X\rvert$ and $y\in h^{-1}(x)$. Since $\kappa(y)$ is finite,
\[\Gal(\kappa(y)/\kappa(x))=\left<\Frob_{\N x}\right>\, ,\]
where $\Frob_{\N x}(a)=a^{\N x}$ for all $a\in \kappa(y)$.
\begin{Def}[Frobenius Element]
    Let $x\in \lvert X\rvert$ be unramified in $Y$, and let $y\in h^{-1}(x)$. We define the Frobenius element at $y/x$ to be
    \[\Fr_{Y/X}[y/x]=\Psi_{y/x}^{-1}(\Frob_{\N x})\in D_{y/x}\subseteq G\, .\]
    It is the unique element $\sigma\in G$ satisfying that for all $a\in \mathcal{O}_{Y,x}$, $\sigma(a)=a^{\N x} \mod m_{y}$.
\end{Def}
\begin{Def}[Frobenius Class]\label{Def:12}
    Let $x\in \lvert X\rvert $ be unramified in $Y$. We define the Frobenius class at $x$ in $Y$ to be
    \[\Fr_{Y/X}(x):=\{\Fr_{Y/X}[y/x]\mid y\in h^{-1}(x)\}\, .\]
\end{Def}
\begin{Lemma}
    Let $x\in \lvert X\rvert$ be unramified in $Y$. Then, $\Fr_{Y/X}(x)$ is a conjugacy class in $G$.
\end{Lemma}
\begin{proof}
    If $\sigma\in \Fr_{Y/X}(x)$ and $\tau\in G$, then
    \[\Psi_{y/x}(\tau\sigma\tau^{-1})=\Psi_{y/x}(\tau)\Psi_{y/x}(\sigma)\Psi_{y/x}(\tau)^{-1}=\Psi_{y/x}(\sigma)=\Fr_{\N x}\, .\]
    Therefore, $\tau\sigma\tau^{-1}\in \Fr_{Y/X}(x)$.

    If $y,y'\in h^{-1}(x)$, then there exists a $\sigma\in G$ that satisfies
    \[\sigma[m_{y}]=m_{y'}\, .\]
    Therefore, since $\sigma D_{y/x}\sigma^{-1}=D_{y'/x}$ and
    \[\sigma\Fr_{Y/x}[y/x]\sigma^{-1}\in \Fr_{Y/X}(x)\, ,\]
    it follows that $\sigma\Fr_{Y/x}[y/x]\sigma^{-1}=\Fr_{Y/x}[y'/x]$.
\end{proof}
\begin{Prop}\label{Prop:1}
    The field $K$ is a regular extension of $\mathbb{F}_{q}$, of transcendence degree $\td(K)=d<\infty$.
\end{Prop}
\begin{proof}
    Since $K=\Frac(\mathcal{O}_{X}(X))$ is the fraction field of a finitely generated $\mathbb{F}_{q}$-algebra, the transcendence degree is finite.
    Since $X_{\overline{\mathbb{F}_{q}}}$ is irreducible, $K\otimes_{\mathbb{F}_{q}} \overline{\mathbb{F}_{q}}=\kappa(X_{\overline{\mathbb{F}_{q}}})$ is a field, meaning that $K$ is a regular extension of $\mathbb{F}_{q}$. Lastly, the equality $\td(K)=d$ is shown in ~\cite[Chapter 2, Proposition 5.19]{QingLiu}.
\end{proof}
\begin{Not}
    For every $r\in \mathbb{N}$, denote by $p_{r,X}$ the projection $p_{r,X}:X_{\mathbb{F}_{q^{r}}}\to X$.
\end{Not}
\begin{Prop}\label{Cor:3}
    Let $j,r\in \mathbb{N}$ be such that $j|r$, and let $x\in \lvert X\rvert_{j}$. Then, under the morphism $p_{r,X}:X_{\mathbb{F}_{q^{r}}}\to X$,
    \[f_{x}=\frac{r}{j},\ e_{x}=1,\ g_{x}=j\, .\]
\end{Prop}
\begin{proof}
    Let $x\in X$ be closed, and let $x\in U\subseteq X$ be an affine scheme neighbourhood of $x$. Take $A=\mathcal{O}_{X}(U)$.
    Since $p_{r,X}^{-1}[U]\cong \Spec A\mathbb{F}_{q^{r}}$, for every $y\in p_{r,X}^{-1}(x)$, $\kappa(y)\cong \mathbb{F}_{q^{r}}$. Thus,
    \[f_{x}=[\mathbb{F}_{q^{r}}:\mathbb{F}_{q^{j}}]=\dfrac{r}{j}\, .\]
    Additionally, from Proposition~\ref{Prop:5}, $x$ is unramified, meaning that $e_{x}=1$. Lastly, applying Corollary~\ref{lemma:7},
    \[\lvert p_{r,X}^{-1}(x)\rvert=\dfrac{r}{e_{x}f_{x}}=j\, .\]
\end{proof}
\begin{Def}
Let $\mathcal{C} \subseteq G$ be a conjugacy class, and let $\sigma\in \mathcal{C}$. We define
\[P_{Y/X}(\mathcal{C})=\{x\in X\mid x \text{ is unramified in } Y,\ \Fr_{Y/X}(x)=\mathcal{C}\}\, ,\]
and
\[P_{Y/X}(\sigma):=\{y\in Y\mid h(y) \text{ is unramified in }Y,\ \Fr_{Y/X}[y/x]=\sigma\}\, .\]
For every $1\leq k\in \mathbb{N}$, we define
\[P_{Y/X}(k,\mathcal{C})=\{x\in \lvert X\rvert_{k} \mid x \text{ is unramified in } Y,\ \Fr_{Y/X}(x)=\mathcal{C}\}\, ,\]
and denote
\[\pi_{Y/X}(k,\mathcal{C})=\lvert P_{Y/X}(k,\mathcal{C}) \rvert\, .\]
\end{Def}
\begin{remark}
Since $\sigma \in \mathcal C$, we have
    \[h(P_{Y/X}(\sigma))=P_{Y/X}(\mathcal C)\, .\]
\end{remark}
\begin{Def}[Degree of an Affine Variety]\label{Def:2.5}
    Let $E/\mathbb{F}_{q}$ be an algebraic extension, $Z$ be an affine $E$-variety, and let $h\colon Z\to \mathbb{A}^{m}_{E}$ be the associated closed immersion into an affine space. The degree of $Z$ is defined as:
    \[\delta_{Z,h}:=\max\{\lvert Z(\overline{\mathbb{F}_{q}})\cap L\rvert\mid L\subseteq \overline{\mathbb{F}_{q}}^{m} \text{ is an affine linear subspace},\ \codim(L)=\dim(Z),\ \lvert Z(\overline{\mathbb{F}_{q}})\cap L\rvert<\infty\}\, .\]
    In most cases we will omit the $h$, and write $\delta_{Z}$ instead.
\end{Def}
\begin{remark}\label{rem:degree}
    An immediate consequence from this definition is that if $Z$ is an affine $\mathbb{F}_{q^{n}}$-variety, then $\delta_{Z}=\delta_{Z_{{\mathbb{F}_{q^{k}}}}}$ for all $k\in \mathbb{N}$ divisible by $n$.
\end{remark}
\section{Explicit Chebotarev's Density Theorem}\label{Section:3} We use Setting~\ref{Setting:first} and Notation~\ref{Setting}, and fix a conjugacy class $\mathcal{C}\subseteq G$.

The general idea of the proof is to find a variety $Z$ for which $\pi_{Y/X}(k,\mathcal{C})$ can be approximated by $\lvert Z(\mathbb{F}_{q})\rvert$. Roughly speaking, this is done by taking some $\sigma\in \mathcal{C}$, taking $Y'$ to be an irreducible component of $Y_{\mathbb{F}_{q^{nk\ord(\sigma)}}}$, and then taking $Z$ to be the normalisation of $X$ in $(L\mathbb{F}_{q^{nk\ord(\sigma)}})^{\sigma}$. We get the following diagram:
\[\begin{tikzcd}
	& {Y'} \\
	Y && {Z} \\
	& X
	\arrow["{p_{nk\ord(\sigma),Y}}"', from=1-2, to=2-1]
	\arrow["{f}", from=1-2, to=2-3]
	\arrow["h"', from=2-1, to=3-2]
	\arrow["g", from=2-3, to=3-2]
\end{tikzcd}\, .\]
\begin{Lemma}\label{lemma:8}
    Let $k\in \mathbb{N}$ be such that for all $\sigma\in \mathcal{C}$, $\sigma\vert_{\mathbb{F}_{q^{n}}}=\Frob_{q}^{k}\vert_{\mathbb{F}_{q^{n}}}$. Let $n' \in \mathbb{N}$ be some multiple of $n$ and take $Y'$ to be an irreducible component of ${Y}\times_{\mathbb{F}_{q}}\mathbb{F}_{q^{n'}}$.
    Denote $L'=\kappa(Y')$. 
    Then,
    \begin{enumerate}[label=\arabic*.]
        \item The map $gf\colon Y'\to X$ is generically Galois, and $\ScDeg(Y')=n'$.
        \item For each $\sigma \in \mathcal{C}$ there exists a unique $\sigma'\in \Gal(L'/K)$ that satisfies \[\sigma'\vert_{L}=\sigma,\ \sigma'\vert_{\mathbb{F}_{q^{n'}}}=\Frob^{k}_{q}\vert_{\mathbb{F}_{q^{n'}}}\, .\]
        \item For each $\sigma \in \mathcal{C}$, $\ord(\sigma')=\lcm(\ord(\sigma), [\mathbb{F}_{q^{n'}}:\mathbb{F}_{q^{n'}}\cap\mathbb{F}_{q^{k}}])$.
        \item 
        $\mathcal{C}'=\{\sigma'\mid\sigma \in \mathcal{C}\}\subseteq \Gal(L'/K)$ is a conjugacy class.
        \item $P_{Y/X}(k, \mathcal{C})=P_{Y'/X}(k, \mathcal{C}')$.
    \end{enumerate}  
\end{Lemma}
\begin{proof}
    Since $Y'$ is an irreducible component of $Y\times_{\mathbb{F}_{q}} \mathbb{F}_{q^{n'}}$, $L'=L\mathbb{F}_{q^{n'}}$, and therefore $L'/K$ is Galois. Additionally, since $L$ is linearly disjoint from $\mathbb{F}_{q^{n'}}$ over $\mathbb{F}_{q^{n}}$, and $n|n'$, we conclude that \[\ScDeg(Y')=[L\mathbb{F}_{q^{n'}}:L][L:\mathbb{F}_{q}]=\dfrac{n'}{n}n=n'\, .\]
    The existence of $\sigma'$ in condition $2$ follows from the linear disjointness of $L$ from $\mathbb{F}_{q^{n'}}$ over $\mathbb{F}_{q^{n}}$.
    
    We now prove $3$. Suppose $\sigma\in \mathcal{C}$. Then, for all $l\in \mathbb{N}$,
    \[\sigma'^{l}=\id\iff \sigma^{l}=\id, \Frob_{q}^{kl}\vert_{\mathbb{F}_{q^{n'}}}=\id\, .\]
    In other words, $\ord(\sigma')\mid l$ if and only if $\ord(\sigma)\mid l$, and $\frac{n'}{\gcd(n',k)}=\ord(\Frob_{q}^{k}\vert_{\mathbb{F}_{q^{n'}}})\mid l$. This means that  \[\ord(\sigma')=\lcm(\ord(\sigma),\frac{n'}{\gcd(n',k)})=\lcm(\ord(\sigma),[\mathbb{F}_{q^{n'}}:\mathbb{F}_{q^{n'}}\cap \mathbb{F}_{q^{k}}])\, .\]
    For $4$, let $\sigma,\tau\in \mathcal{C}$, and take a $\varphi\in G$ that satisfies $\sigma=\varphi\tau\varphi^{-1}$. Then, there is an $s\in \mathbb{N}$ for which
    \[\varphi\vert_{\mathbb{F}_{q^{n}}}=\Frob_{q}^{s}\vert_{\mathbb{F}_{q^{n'}}}\, .\]
    In a similar fashion as $2$, from linear disjointness of $L$ from $\mathbb{F}_{q^{n'}}$ over $\mathbb{F}_{q^{n}}$, there exists a $\varphi'\in \Gal(L'/K)$ that satisfies \[\varphi'\vert_{L}=\varphi,\ \varphi'\vert_{\mathbb{F}_{q^{n'}}}=\Frob_{q}^{s}\vert_{\mathbb{F}_{q^{n'}}}\, ,\] then $\sigma'=\varphi'\tau'\varphi'^{-1}$.
    
    Lastly, we prove $5$. By Proposition~\ref{Cor:transitive} and Proposition~\ref{Prop:5}, $x\in \lvert X\rvert$ is unramified in $Y$ if and only if $x$ is unramified in $Y'$.
    Let $y'\in (hp_{n',Y})^{-1}(x)$, $y=p_{n',Y}(y)\in h^{-1}(x)$, and take $\sigma=\Fr_{Y/X}[y/x]$. For all $a\in \mathcal{O}_{Y,x}$ and $u\in \mathbb{F}_{q^{n'}}$ we have
    \[\sigma'(au)=\sigma(a)\Frob_{q}^{k}(u)=a^{q^{k}}u^{q^{k}} \mod m_{y'}\, .\]
     Since $\mathcal{O}_{Y',x}=\mathcal{O}_{Y,x}\mathbb{F}_{q^{n'}}$, it follows that for all $b\in \mathcal{O}_{Y'x}$,\[\sigma'(b)=b^{q^{k}}\mod m_{y'}\, ,\]
     meaning that $\sigma'=\Fr_{Y'/X}[y'/x]$.
\end{proof}
\begin{Lemma}\label{lemma:9}
    Let $x\in P_{Y/X}(\mathcal{C})$, $\sigma \in \mathcal{C}$. Then,\\
    \[g_{x} = \dfrac{\lvert G\rvert}{\ord(\sigma)}\,,  \lvert\{y\in h^{-1}(x)\mid \Fr_{Y/X}[y/x]=\sigma\}\rvert=\dfrac{\lvert G\rvert}{\lvert \mathcal{C}\rvert \ord(\sigma)}\, .\]
    In particular, if $\mathcal{F}\subseteq P_{Y/X}(\mathcal{C})$ is finite, then
    \[\lvert\mathcal{F}\rvert=\dfrac{\lvert\mathcal{C}\rvert}{\lvert G\rvert}\ord(\sigma)\lvert\{y\in h^{-1}(\mathcal{F})\mid \Fr_{Y/X}[y/h(y)]=\sigma\}\rvert\, .\]
\end{Lemma}
\begin{proof}
    By Corollary~\ref{lemma:7}, we have
    \[\lvert G\rvert=e_{x}f_{x}g_{x}=f_{x}g_{x}\, .\]
    Additionally, since $x$ is unramified and $\sigma\in \mathcal{C}=\Fr_{Y/X}(x)$:
    \[e_{x}=1,\ f_{x}=\ord(\sigma)\, .\]
    Therefore, \[g_{x}=\dfrac{\lvert G\rvert}{\ord(\sigma)}\, .\]
    For the second part, notice that for all $\tau\in \mathcal{C}$, we can choose a  $\mu\in G$ that satisfies $\tau=\mu\sigma\mu^{-1}$. Then, if $y\in h^{-1}(x)$ satisfies $\Fr_{Y/X}[y/h(y)]=\sigma$, there is a unique $y'\in h^{-1}(x)$ that satisfies $\mu\sigma[m_{y}]=m_{y'}$. Therefore, this $y'$ satisfies:
    \[\Fr_{Y/X}[y'/x]=\mu\Fr_{Y/X}[y/x]\mu^{-1}=\mu\sigma\mu^{-1}=\tau\, .\]
    Thus, since the map $y\leftrightarrow y'$ is bijective,
    \[\lvert \{y\in h^{-1}(x)\mid \Fr_{Y/X}[y/h(y)]=\sigma\}\rvert=\lvert \{y\in h^{-1}(x)\mid \Fr_{Y/X}[y/h(y)]=\tau\}\rvert\, .\]
    Therefore,
    \[
        \dfrac{\lvert G\rvert}{\ord(\sigma)}=\lvert h^{-1}(x)\rvert=\lvert\biguplus_{\tau \in \mathcal{C}} \{y\in h^{-1}(x)\mid \Fr_{Y/X}[y/x]=\tau\}\rvert
        =\lvert \mathcal{C}\rvert \lvert\{y\in h^{-1}(x)\mid \Fr_{Y/X}[y/x]=\sigma\}\rvert\, .
    \]
\end{proof}
\begin{Lemma}\label{lemma:10}
    Let $Z$ be an $\mathbb{F}_{q}$-variety, with morphisms 
    \[\begin{tikzcd}
	{Y} & {Z} & {X}
	\arrow["{f}", from=1-1, to=1-2]
	\arrow["h", curve={height=-24pt}, from=1-1, to=1-3]
	\arrow["g", from=1-2, to=1-3]
\end{tikzcd}\, .\]
    Denote $r=\ScDeg(Z)$ and $F=\kappa(Z)$. Let $r|k$, and let $y\in Y$ be a closed point such that $h(y)$ is unramified in $Y$, and $\deg(h(y))=k$. Then:
    \[\Fr_{Y/X}[y/h(y)]\in \Gal(L/F) \text{ if and only if } \deg(f(y))=k/r\, .\]
    In this case, $\Fr_{Y/X}[y/h(y)]=\Fr_{Y/Z}[y/f(y)]$.
\end{Lemma}
\begin{proof}
    By definition, for all $b\in \kappa(y)$,
    \[\Psi_{y/h(y)}(\Fr_{Y/X}[y/h(y)])(b)=b^{q^{k}}\]
    and
    \[\Psi_{y/f(y)}(\Fr_{Y/X}[y/f(y)])(b)=b^{q^{r\deg(f(y))}}\, .\]
    Therefore, if $\deg(f(y))=\frac{k}{r}$, then
    \[\Psi_{y/h(y)}(\Fr_{Y/X}[y/h(y)])=\Psi_{y/f(y)}(\Fr_{Y/Z}[y/f(y)])=\Psi_{y/h(y)}(\Fr_{Y/Z}[y/f(y)])\, ,\]
    and since $h(y)$ is unramified, $\Psi_{y/h(y)}$ is injective, and
    \[\Fr_{Y/X}[y/h(y)]=\Fr_{Y/Z}[y/f(y)]\, .\]
    Consequently, $\Fr_{Y/X}[y/h(y)]\in \Gal(L/F)$.
    
    Conversely,
    if, $\Fr_{Y/X}[y/h(y)]\in \Gal(L/F)$, then for all $a\in \kappa(f(y))$,
    \[\Fr_{Y/X}[y/h(y)](a)=a\, ,\]
    and so $a=a^{q^{k}}$.
    Therefore,
    \[\mathbb{F}_{q^{r\deg(f(y))}}=\kappa(f(y))\subseteq\mathbb{F}_{q^{k}}\, .\]
    However,
    \[\mathbb{F}_{q^{k}}=\kappa(h(y))\subseteq \kappa(f(y))=\mathbb{F}_{q^{r\deg(f(y))}}\, ,\]
    so $\deg(f(y))=\frac{k}{r}$.
\end{proof}
In the same notations, we can now express $\pi_{Y/X}(k,\mathcal{C})$ using $P_{Y/Z}(\mathcal{C}')$:
\begin{Lemma}\label{lemma:11}
    With the same notations as in Lemma~\ref{lemma:10}, let $\sigma \in \Gal(L/F)$, and denote by $\mathcal{C}, \tilde{\mathcal{C}}$ its conjugacy classes in $G, \Gal(L/F)$, respectively. Then,\\
    \[\pi_{Y/X}(k,\mathcal{C}) = \dfrac{\lvert \mathcal{C}\rvert}{\lvert \tilde{\mathcal{C}}\rvert [F:K]} \lvert\{z\in P_{Y/Z}(k/r,\tilde{\mathcal{C}})\mid g(z)\text{ is unramified},\ \deg_{X}(g(z))>k/2\}\rvert\, .\]
\end{Lemma}
\begin{proof}
Let $z\in P_{Y/Z}(k/r,\tilde{\mathcal{C}})$ such that $\deg(g(z))>\frac{k}{2}$ Then,
\[\mathbb{F}_{q^{\deg(g(z))}}= \kappa(g(z))\subseteq \kappa(z)= \mathbb{F}_{q^{r\frac{k}{r}}}=\mathbb{F}_{q^{k}}\, .\]
Therefore, $\deg(g(z))=k$, as $\deg(g(z))|k$ and $\deg(g(z))>\frac{k}{2}$.
Now, denote
\[\mathcal{F}'=\{z\in P_{Y/Z}(k/r,\tilde{\mathcal{C}})\mid g(z)\text{ is unramified in } Z,\deg(g(z))=k\}\, .\]
Applying Lemma~\ref{lemma:9} gives
\[\lvert \mathcal{F}'\rvert=\dfrac{\lvert \tilde{\mathcal{C}}\rvert}{[L:F]}\ord(\sigma)\lvert \{y\in f^{-1}(\mathcal{F}')\mid\Fr_{Y/X}[y/h(y)]=\sigma\}\rvert\, .\]
From Lemma~\ref{lemma:10} we get
    \[\{y\in f^{-1}(\mathcal{F}')\mid \Fr_{Y/X}[y/h(y)]=\sigma\}
    =\{y\in P_{Y/Z}(\sigma)\mid \deg(h(y))=k\}\, .\]
Finally, apply Lemma~\ref{lemma:9} with $\mathcal{F}=P_{Y/X}(k,\mathcal{C})$ to obtain
\begin{align*}
    \pi_{Y/X}(k,\mathcal{C})=\lvert \mathcal{F}\rvert&=\dfrac{\lvert \mathcal{C}\rvert}{[L:K]}\ord(\sigma)\lvert \{y\in h^{-1}(\mathcal{F})\mid h(y)\text{ is unramified, } \Fr_{Y/X}[y/h(y)]=\sigma\}\rvert\\
    &=\dfrac{\lvert \mathcal{C}\rvert}{[L:K]}\ord(\sigma)\lvert \{y\in P_{Y/Z}(\sigma)\mid \deg(h(y))=k\}\rvert\\
    &=\dfrac{\lvert \mathcal{C}\rvert}{[L:K]}\ord(\sigma)\dfrac{\lvert \mathcal{F}'\rvert[L:F]}{\ord(\sigma)\lvert \tilde{\mathcal{C}}\rvert}=\dfrac{\lvert \mathcal{C}\rvert \lvert \mathcal{F}'\rvert}{\lvert \tilde{\mathcal{C}}\rvert[F:K]}\, .
\end{align*}
\end{proof}
The following lemma approximates $\pi_{X_{\mathbb{F}_{q^{n}}}/X}(1,\mathcal{C})$.
\begin{Lemma}\label{lemma:15}
    Assume that $q>2(d+1)\delta_{X}^{2}$. If $\sigma\in \Gal(K\mathbb{F}_{q^{n}}/K)$ satisfies that $\sigma\vert_{\mathbb{F}_{q^{n}}}=\Frob_{q}\vert_{\mathbb{F}_{q^{n}}}$, then
    \[\lvert\pi_{X_{\mathbb{F}_{q^{n}}}/X}(1,\{\sigma\})-q^{d}\rvert\leq (\delta_{X}-1)(\delta_{X}-2)q^{d-\frac{1}{2}}+5\delta_{X}^{\frac{13}{3}}q^{d-1}\, .\]
\end{Lemma}
\begin{proof}
    By Proposition~\ref{Prop:5}, every closed point $x\in X$ is unramified in $X_{\mathbb{F}_{q^{n}}}$. Therefore,
    \[P_{X_{\mathbb{F}_{q^{n}}}/X}(1,\{\sigma\})=\lvert X\rvert_{1}\, .\]
    Thus, \cite[Theorem 7.1]{cafurematera} implies
    \[\lvert\pi_{X_{\mathbb{F}_{q^{n}}}/X}(1,\{\sigma\})-q^{d}\rvert\leq (\delta_{X}-1)(\delta_{X}-2)q^{d-\frac{1}{2}}+5\delta_{X}^{\frac{13}{3}}q^{d-1}\, ,\]
    as needed.
\end{proof}
We also have the following bound:
\begin{Lemma}\label{lemma:15'}
    Let $k\in \mathbb{N}$, and let $Z$ be a geometrically irreducible $\mathbb{F}_{q}$-variety of dimension $d$, that is also an affine scheme. Denote $K'=\kappa(Z)$. Assume there exists a closed immersion
    \[\iota\colon Z_{\overline{\mathbb{F}_{q}}}\hookrightarrow \mathbb{A}^{v}_{\overline{\mathbb{F}_{q}}}\, ,\]
    for some $v\in \mathbb{N}$, such that
     $\iota(Z_{\overline{\mathbb{F}_{q}}})$ is a variety of degree $\delta_{Y_{\overline{\mathbb{F}_{q}}}}$, defined by the vanishing of $l$ polynomials, of degrees bounded by $D$.
    Then,
        For all $n'\in \mathbb{N}$, if $\sigma\in \Gal(K'\mathbb{F}_{q^{n'}}/K')$ satisfies that $\sigma\vert_{\mathbb{F}_{q^{n'}}}=\Frob_{q}\vert_{\mathbb{F}_{q^{n'}}}$, then
        \[\lvert\pi_{Z_{\mathbb{F}_{q^{n'}}}/Z}(1,\{\sigma\})-q^{d}\rvert\leq (\delta_{Y_{\overline{\mathbb{F}_{q}}}}-1)(\delta_{Y_{\overline{\mathbb{F}_{q}}}}-2)q^{d-\frac{1}{2}}+6\cdot 2^{l}(lD+3)^{v+1}q^{d-1}\, .\]
\end{Lemma}
\begin{proof}
    As in the proof of Lemma~\ref{lemma:15},
    \[P_{Z_{\mathbb{F}_{q^{n'}}}/Z}(1,\{\sigma\})=\lvert Z\rvert_{1}\, .\]
    Thus, since the proof of \cite[Theorem 11.1]{GL} depends only on constants of $Z_{\overline{\mathbb{F}_{q}}}$, combining it with $(1)$ and \cite[Remark 11.3]{GL} implies
    \[\lvert\pi_{Z_{\mathbb{F}_{q^{n'}}}/Z}(1,\{\sigma\})-q^{d}\rvert \leq (\delta_{Y_{\overline{\mathbb{F}_{q}}}}-1)(\delta_{Y_{\overline{\mathbb{F}_{q}}}}-2)q^{d-\frac{1}{2}}+6\cdot 2^{l}(lD+3)^{v+1}q^{d-1}\, ,\]
    as needed.
\end{proof}
\begin{Lemma}\label{lemma:14}
    Suppose that $q>2(d+1)\delta_{X}^{2}$. Let $k\in \mathbb{N}$, and let
    \[\begin{tikzcd}
	{Z} & {X}
    \arrow["{f}", from=1-1, to=1-2]
    \end{tikzcd}\]
    be a finite morphism of normal $\mathbb{F}_{q}$-varieties.
    
    Denote $K'=\kappa(Z)$. Assume that $K'/K$ is finite and separable, $[K':K]=km$ and $\mathbb{F}_{q^{k}}\subseteq K'$.
    
    Then,
    \begin{align*}
        \lvert \{z\in Z\mid f(z) \text{ is closed, } k\neq \deg_{X}(f(z))\vert k\}\rvert\leq 2m(\delta_{X}-1)(\delta_{X}-2)&q^{\frac{k}{2}(d-\frac{1}{2})}\\
        +10m\delta_{X}^{\frac{13}{3}}&q^{\frac{k}{2}(d-1)}\\
        +2m&q^{\frac{k}{2}d}
    \end{align*}
\end{Lemma}
\begin{proof}
    Let $j|k$ be some divisor of $k$, and let $x\in X$ be of degree $\deg(x)=j$. First we prove that
    \[\lvert f^{-1}(x)\rvert \leq [K':K\mathbb{F}_{q^{k}}]j\]
    Indeed, factor $f\colon Z\to X$ as
    \[\begin{tikzcd}
	{Z} & {X_{ \mathbb{F}_{q^{k}}}} & X
	\arrow["{\tilde{f}}", from=1-1, to=1-2]
	\arrow["f", curve={height=-28pt}, from=1-1, to=1-3]
	\arrow["{p_{k,X}}", from=1-2, to=1-3]
\end{tikzcd}\]
    By Proposition~\ref{Cor:3},
    \[\lvert p_{k,X}^{-1}(x)\rvert=j\, .\]
    By Lemma~\ref{lemma:3.6}, for every $\tilde{x}\in \lvert X_{ \mathbb{F}_{q^{k}}}\rvert $,
    \[\lvert \tilde{f}^{-1}(\tilde{x})\rvert\leq [K':K\mathbb{F}_{q^{k}}]\, .\]
    Therefore,
    \[\lvert f^{-1}(x)\rvert\leq [K':K\mathbb{F}_{q^{k}}]j\, ,\]
    as wanted.
    Now, since
    \[j=\lvert \{x'\in X_{\mathbb{F}_{q^{j}}}\mid p_{k,X}(x')=x\}\rvert\, ,\]
    we know that
    \[\lvert f^{-1}(x)\rvert \leq [K':K\mathbb{F}_{q^{k}}]\lvert \{x'\in X_{\mathbb{F}_{q^{j}}}\mid p_{k,X}(x')=x\}\rvert\, .\]
    Using this, and the fact that every $x'\in p_{k,X}^{-1}(x)$ is of degree $\deg(x')=1$, we get that

    \begin{align*}
        \lvert \{z\in Z\mid f(z) \text{ is closed, }& k\neq \deg_{X}(f(z))\vert k\}\rvert 
        = \sum_{\substack{j|k \\ j\neq k}} \lvert \{z\in Z\mid f(z) \text{ is closed, } \deg_{X}(f(z))=j\}\rvert \\
        &\qquad\leq m\sum_{\substack{j|k \\ j\neq k}}\lvert \{x'\in X_{\mathbb{F}_{q^{j}}}\mid p_{k,X}(x') \text{ is closed, } \deg_{X}(p_{k,X}(x'))=j\}\rvert\\
        &\qquad\leq m\sum_{\substack{j|k \\ j\neq k}}\lvert X_{\mathbb{F}_{q^{j}}}\rvert_{1}\, .
    \end{align*}
    In \cite[Theorem 7.1]{cafurematera} it is shown that for all $k\neq j|k$,
    \begin{align*}
       \lvert X_{\mathbb{F}_{q^{j}}}\rvert_{1}\leq (\delta_{X_{\mathbb{F}_{q^{j}}}}-1)(\delta_{X_{\mathbb{F}_{q^{j}}}}-2)&q^{j(d-\frac{1}{2})}\\
    +5\delta_{X_{\mathbb{F}_{q^{j}}}}^{\frac{13}{3}}&q^{j(d-1)}\\
    +&q^{jd}\, .
    \end{align*}
    Since for all $1\leq l$,
    \[\sum_{\substack{j|k \\ j\neq k}} q^{jl}\leq 2q^{\frac{kl}{2}}\, ,\]
    and since $\delta_{X_{\mathbb{F}_{q^{j}}}}=\delta_{X}$ (see remark~\ref{rem:degree}), we get that for $d>1$:
    \begin{multline*}
        \lvert \{z\in Z\mid f(x') \text{ is closed, } k\neq \deg_{X}(f(z))\vert k\}\rvert\\
        \leq m\sum_{\substack{j|k \\ j\neq k}}\left[(\delta_{X}-1)(\delta_{X}-2)q^{j(d-\frac{1}{2})} +5\delta_{X}^{\frac{13}{3}}q^{j(d-1)}
        +q^{jd}\right]\\
        \leq 2m(\delta_{X}-1)(\delta_{X}-2)q^{\frac{k}{2}(d-\frac{1}{2})}+10m\delta_{X}^{\frac{13}{3}}q^{\frac{k}{2}(d-1)}+2mq^{\frac{k}{2}d}\, ,
    \end{multline*}
    as needed.
    For $d=1$, the desired bound can be achieved by applying \cite[Theorem 5.5.2]{friedjarden} to bound $\lvert X_{\mathbb{F}_{q^{j}}}\rvert_{1}$ and by using $\sum_{\substack{j|k\\ j\neq k}}(q^{j}+1)\leq 2q^{\frac{k}{2}}$ (which holds since $q>2(d+1)\delta_{X}^{2}\geq 4$).
\end{proof}
We are now ready to prove Theorem~\ref{Thm:1.2}:
\begin{Thm*}\label{Thm:1.1}
    Assume that $X$ is locally Noetherian, and that $q>2(d+1)\delta_{X}^{2}$.
    Let $r\in \mathbb{N}$ be such that for all $\sigma\in \mathcal{C}$, $\sigma\vert_{\mathbb{F}_{q^{n}}}=\Frob_{q}^{r}\vert_{\mathbb{F}_{q^{n}}}$, and let $k\in \mathbb{N}$ be a natural number.

    If $k=r\mod n$, then:
    \begin{align*}
    \bigg|\dfrac{\pi_{Y/X}(k,\mathcal{C})}{q^{kd}/k}-\frac{\lvert \mathcal{C}\rvert}{m}\bigg|&\leq 2\lvert\mathcal{C}\rvert\left[(\delta_{X}-1)(\delta_{X}-2)q^{-\frac{k}{2}(d+\frac{1}{2})}+5\delta_{X}^{\frac{13}{3}}q^{-\frac{k}{2}(d+1)}+q^{-\frac{k}{2}d}\right]\\
        &\quad+\dfrac{\lvert\mathcal{C}\rvert}{m}\left[(\delta_{Y_{\overline{\mathbb{F}_{q^{n}}}}}-1)(\delta_{Y_{\overline{\mathbb{F}_{q^{n}}}}}-2)q^{-\frac{k}{2}}+6\cdot 2^{l}(lD+3)^{v+1}q^{-k}\right]+\lvert \mathcal{C}\rvert R_{k}kq^{-kd}\, ,
    \end{align*}
    Where $R_{k}=\lvert\{x\in X\mid x \text{ is closed and ramified in } Y,\ \deg_{X}(x)=k\}\rvert$.
    
    Otherwise, $\pi_{Y/X}(k,\mathcal{C})=0$.
\end{Thm*}
\begin{proof}
    If $k\neq r \mod n$, and $x\in P_{Y/X}(k, \mathcal{C})$, then for all $y\in h^{-1}(x)$, we get:
    \[\Frob_{q}^{k}\vert_{\mathbb{F}_{q^{n}}}=\Fr_{Y/X}[y/x]\vert_{\mathbb{F}_{q^{n}}}=\Frob_{q}^{r}\vert_{\mathbb{F}_{q^{n}}}\, ,\]
    which is a contradiction.

    Assume that $k=r \mod n$. Choose some $\sigma \in \mathcal{C}$, and denote $n'=nk\ord(\sigma)=nko$. Take
    $Y'$ to be an irreducible component of ${Y}\times_{\mathbb{F}_{q}}\mathbb{F}_{q^{n'}}$,
    and denote $L'=\kappa(Y')$.
    
    Then, from Lemma~\ref{lemma:8}, $L'/K$ is a Galois extension, and there exists some $\sigma'\in \Gal(L'/K)$ such that:
    \[\sigma'\vert_{L}=\sigma,\text{ } \sigma'\vert_{\mathbb{F}_{q^{n'}}}=\Frob_{q}^{k}\vert_{\mathbb{F}_{q^{n'}}}\, .\]
    Denote
    \[K'=L'^{\sigma'}=\{\alpha\in L'\mid \sigma'(\alpha)=\alpha\}\, ,\]
    and let
    \[g\colon Z\to X\]
    be the normalisation of $X$ in $K'$. $X$ is an affine scheme, and therefore $Z$ is an affine scheme as well. Since $X$ is locally Noetherian, $hp_{n',Y}\colon Y'\to X$ is finite and $g$ is integral, we conclude that $f$ is finite.
    
    Any normalisation of $Z$ in $L'$ will be a normalisation of $X$ in $L'$, and therefore there exists an integral morphism $f\colon Y'\to Z$ such that the following diagram commutes:
    \[\begin{tikzcd}
	& {Y'} \\
	Y && {Z} \\
	& X
	\arrow["p_{n',Y}", from=1-2, to=2-1]
	\arrow["{f}"', from=1-2, to=2-3]
	\arrow["h", from=2-1, to=3-2]
	\arrow["g"', from=2-3, to=3-2]
\end{tikzcd}\]
   Since $Y'$ is covered by spectra of finitely generated $\mathbb{F}_{q}$-algebras, $f$ is not only integral, but also finite.

    We now show that $[K':K]=km$ and $L'=K'\mathbb{F}_{q^{n'}}$.
    Notice that \[K'\cap \overline{\mathbb{F}_{q}}=K'\cap L'\cap \overline{\mathbb{F}_{q}}=K'\cap \mathbb{F}_{q^{n'}}\, .\]
    Therefore, if $a\in K'\cap \overline{\mathbb{F}_{q}}$, then $x=\sigma'(x)=x^{q^{k}}$ and so $x\in \mathbb{F}_{q^{k}}$. Thus,
    \[K'\cap \overline{\mathbb{F}_{q}}=\mathbb{F}_{q^{k}}\, .\]
    Additionally, from Lemma~\ref{lemma:8}, \[[L':K']=\ord(\sigma')=\lcm(o,[\mathbb{F}_{q^{n'}}:\mathbb{F}_{q^{k}}])=[\mathbb{F}_{q^{n'}}:\mathbb{F}_{q^{k}}]=no=[K\mathbb{F}_{q^{n'}}:K\mathbb{F}_{q^{k}}]\, .\]
    Therefore,
    \[[K':K]=\frac{[L':L][L:K]}{[L':K']}=\frac{konm}{no}=km\, ,\]
    and since
    \[[K'\mathbb{F}_{q^{n'}}:K']=[\mathbb{F}_{q^{n'}}:\mathbb{F}_{q^{k}}]=no=[L':K']\, ,\]
    we see that $L'=K'\mathbb{F}_{q^{n'}}$.
    
    Moreover, since $Z\times_{\mathbb{F}_{q^{k}}} \mathbb{F}_{q^{n'}}$ and $Y'$ are normalisations of $Z$ in $L'$, we see that
    \[Y'\cong Z\times_{\mathbb{F}_{q^{k}}} \mathbb{F}_{q^{n'}}\]
    (\cite[\href{https://stacks.math.columbia.edu/tag/037Z}{Tag 037Z}]{stacks-project} gives a proof of this fact).
    
    Let $\mathcal{C}'$ be the conjugacy class of $\sigma'\in \Gal(L'/K)$. From Lemma~\ref{lemma:8}, $P_{Y/X}(k,\mathcal{C})=P_{Y'/X}(k,\mathcal{C}')$.
    Using Lemma~\ref{lemma:11} for $L'/K$, we obtain:
    \[\pi_{Y'/X}(k,\mathcal{C}')=\frac{\lvert \mathcal{C}'\rvert}{[K':K]}\lvert\{z\in P_{Y'/Z}(1,\{\sigma'\})\mid g(z) \text{ is unramified, }\deg_{X}(g(z))>\frac{k}{2}\}\rvert\, ,\]
    since the conjugacy class of $\sigma'$ in $\Gal(L'/K')$ is $\{\sigma'\}$.

    Apply Lemma~\ref{lemma:14} to get:
    \begin{align*}
        \lvert \pi_{Y'/X}(k,\mathcal{C}')-\frac{\lvert \mathcal{C}'\rvert}{[K':K]}\pi_{Y'/Z}(1,\{\sigma'\})\rvert&\leq \frac{\lvert \mathcal{C}'\rvert }{[K':K]}\Bigl[\lvert\{z\in Z\mid z, g(z)\text{ are closed, }k\neq \deg_{X}(g(z))|k\}\rvert\\
        &\hphantom{\leq \frac{\lvert \mathcal{C}'\rvert }{[K':K]}}+\lvert\{z\in Z\mid g(z) \text{ is ramified in } Y', \deg_{X}(g(z))=k\}\rvert\Bigr]\\
        &\leq \frac{2\lvert\mathcal{C}'\rvert}{k}\left[(\delta_{X}-1)(\delta_{X}-2)q^{\frac{k}{2}(d-\frac{1}{2})}\right.\\
        &\hphantom{\leq\frac{2\lvert\mathcal{C}'\rvert}{k}}\quad\left.+5\delta_{X}^{\frac{13}{3}}q^{\frac{k}{2}(d-1)}+q^{\frac{k}{2}d}\right]\\
        &\qquad+\dfrac{\lvert \mathcal{C}'\rvert}{[K':K]}[K':K]R_{k}\, .
    \end{align*}
    Since $Y'\cong Z\times_{\mathbb{F}_{q^{k}}} \mathbb{F}_{q^{n'}}$ and $Z$ is an affine scheme, Lemma~\ref{lemma:15'} implies
    \[\lvert\pi_{Y'/Z}(1,\{\sigma\})-q^{kd}\rvert\leq (\delta_{Y_{\overline{\mathbb{F}_{q^{n}}}}}-1)(\delta_{Y_{\overline{\mathbb{F}_{q^{n}}}}}-2)q^{k(d-\frac{1}{2})}++6\cdot 2^{l}(lD+3)^{v+1}q^{k(d-1)}\, ,\]
    which gives us:
    \begin{align*}
        \bigg|\dfrac{\pi_{Y/X}(k,\mathcal{C})}{q^{kd}/k}-\frac{\lvert \mathcal{C}\rvert}{m}\bigg|&=\dfrac{k}{q^{kd}}\lvert\pi_{Y/X}(k,\mathcal{C})-\frac{\lvert \mathcal{C}\rvert}{[K':K]}q^{kd}\rvert\\
        &=\dfrac{k}{q^{kd}}\lvert\pi_{Y'/X}(k,\mathcal{C}')-\frac{\lvert \mathcal{C}'\rvert}{[K':K]}q^{kd}\rvert\\
        &\leq \dfrac{k}{q^{kd}}\left[\lvert \pi_{Y'/X}(k,\mathcal{C}')-\dfrac{\lvert \mathcal{C}'\rvert}{[K':K]}\pi_{Y'/Z}(1,\{\sigma'\})\rvert\right.\\
        &\hphantom{\leq \dfrac{k}{q^{kd}}}\left.+\dfrac{\lvert \mathcal{C}'\rvert}{[K':K]}\lvert \pi_{Y'/Z}(1,\{\sigma'\})-q^{kd}\rvert\right]\\
        &\leq \dfrac{k}{q^{kd}}\Bigg[\frac{2\lvert\mathcal{C}'\rvert}{k}\left[(\delta_{X}-1)(\delta_{X}-2)q^{\frac{k}{2}(d-\frac{1}{2})}
        +5\delta_{X}^{\frac{13}{3}}q^{\frac{k}{2}(d-1)}+q^{\frac{k}{2}d}\right]+\lvert \mathcal{C}'\rvert R_{k}\\
        &\hphantom{\hphantom{\leq \dfrac{k}{q^{kd}}}}+\dfrac{\lvert\mathcal{C}'\rvert}{[K':K]}\left[(\delta_{Y_{\overline{\mathbb{F}_{q^{n}}}}}-1)(\delta_{Y_{\overline{\mathbb{F}_{q^{n}}}}}-2)q^{k(d-\frac{1}{2})}+6\cdot 2^{l}(lD+3)^{v+1}q^{k(d-1)}\right]\Bigg]\, .
    \end{align*}
    Therefore,
    \begin{align*}
    \bigg|\dfrac{\pi_{Y/X}(k,\mathcal{C})}{q^{kd}/k}-\frac{\lvert \mathcal{C}\rvert}{m}\bigg|&\leq 2\lvert\mathcal{C}\rvert\left[(\delta_{X}-1)(\delta_{X}-2)q^{\frac{k}{2}(-d-\frac{1}{2})}+5\delta_{X}^{\frac{13}{3}}q^{\frac{k}{2}(-d-1)}+q^{-\frac{k}{2}d}\right]\\
    &\quad+\dfrac{\lvert\mathcal{C}\rvert}{m}\left[(\delta_{Y_{\overline{\mathbb{F}_{q^{n}}}}}-1)(\delta_{Y_{\overline{\mathbb{F}_{q^{n}}}}}-2)q^{-\frac{k}{2}}+6\cdot 2^{l}(lD+3)^{v+1}q^{-k}\right]+k\lvert \mathcal{C}\rvert R_{k}q^{-kd}\, .
    \end{align*}
\end{proof}
A slight alteration to this proof will give us Theorem~\ref{Thm:1.3}:
\begin{Thm*}
    Assume that $X$ is locally Noetherian, and that for all $\sigma\in \mathcal{C}$, $\sigma\vert_{\mathbb{F}_{q^{n}}}=\Frob_{q}\vert_{\mathbb{F}_{q^{n}}}$. Then:
    \begin{align*}
    \bigg|\dfrac{\pi_{Y/X}(1,\mathcal{C})}{q^{d}}-\frac{\lvert \mathcal{C}\rvert}{m}\bigg|&\leq \dfrac{\lvert\mathcal{C}\rvert}{m}\left[(\delta_{Y_{\overline{\mathbb{F}_{q^{n}}}}}-1)(\delta_{Y_{\overline{\mathbb{F}_{q^{n}}}}}-2)q^{-\frac{1}{2}}+6\cdot 2^{l}(lD+3)^{v+1}q^{-1}\right]+\lvert \mathcal{C}\rvert R_{1}q^{-d}\, ,
        \end{align*}
    where $R_{1}=\lvert\{x\in X\mid x \text{ is closed and ramified in } Y,\ \deg_{X}(x)=1\}\rvert$.
\end{Thm*}
\begin{proof}
    Following the same steps as in the proof of Theorem~\ref{Thm:1.2}, one can see that
    \begin{align*}
        \lvert \pi_{Y'/X}(1,\mathcal{C}')-\frac{\lvert \mathcal{C}'\rvert}{[K':K]}\pi_{Y'/Z}(1,\{\sigma'\})\rvert&\leq \frac{\lvert \mathcal{C}'\rvert }{[K':K]}\Bigl[\lvert\{z\in Z\mid z, g(z)\text{ are closed, }\deg_{X}(g(z))\leq\frac{1}{2}\}\rvert\\
        &\hphantom{\leq \frac{\lvert \mathcal{C}'\rvert }{[K':K]}}+\lvert\{z\in Z\mid g(z) \text{ is ramified in } Y', \deg_{X}(g(z))=1\}\rvert\Bigr]\\
        &\leq \lvert\mathcal{C}'\rvert R_{1}\, .
    \end{align*}
    Continuing in the same way as in the proof of Theorem~\ref{Thm:1.2} will give the mentioned inequality.
\end{proof}
We now prove a similar bound for The Chebotarev Density Theorem for function fields of one variable.
\section{Applications}\label{Section:5}
\subsection{Polynomial Version of Theorem~\ref{Thm:1.2}}
We provide a version of Theorem~\ref{Thm:1.2} in the case where $Y$ is defined via a polynomial equation.
\begin{Setting}\label{Setting:app}
Let $f\in \mathbb{F}_{q}[x_{1},\ldots,x_{d}][t]$ be a polynomial of degree $N=\deg(f)>1$ that decomposes into monic, separable polynomials in $t$. Let $L$ be the splitting field of $f$ over $\mathbb{F}_{q}(x_{1},\ldots,x_{d})$.
Take \[X=\mathbb{A}_{\mathbb{F}_{q}}^{d}=\Spec \mathbb{F}_{q}[x_{1},\ldots,x_{d}]\, ,\] and let $Y$ be the normalisation of $X$ in $L$. We additionally use Setting~\ref{notation1}.
\end{Setting}
Denote $D=(8N\log(N)+n)Nm(2m-1)-3,\ \delta=DNm$.
We make use of the following proposition about varieties of unramified points.
\begin{Prop}\label{Prop:4.1}
    Let $U_{Y}\subseteq Y, U_{X}\subseteq X$ be the open subsets of all unramified points in $Y$ and $X$ respectively. Then, there exist open subsets $\tilde{U}_{X}\subseteq U_{X}$, $\tilde{U}_{Y}\subseteq U_{Y}$ such that $\tilde{U}_{Y}=h^{-1}[\tilde{U}_{X}]$, and: 
    \begin{enumerate}[label=\arabic*.]
        \item There exists an embedding of ${\tilde{U}_{Y}{}}_{\overline{\mathbb{F}_{q}}}$ into the affine space $\mathbb{A}^{d+2}_{\overline{\mathbb{F}_{q}}}$, such that ${\tilde{U}_{Y}{}}_{\overline{\mathbb{F}_{q}}}$ is defined by the vanishing of $2$ polynomials of degrees bounded by $D$, and $\delta_{{\tilde{U}_{Y}{}}_{\overline{\mathbb{F}_{q}}}}\leq \delta$.
        \item \[\lvert U_{X}\setminus \tilde{U}_{X}\rvert_{1}\leq (D-1)q^{d-1}\, .\]
    \end{enumerate}
\end{Prop}
\begin{proof}
    This proposition is not the main focus of this section, so we prove it in the \hyperref[Appendix]{Appendix}.
\end{proof}
\begin{remark}
    Since $h\colon Y\to X$ is generically Galois, $U_{Y}=h^{-1}[U_{X}]$.
\end{remark}
\begin{Prop}\label{Prop:PolyCheb}
    Let $\mathcal{C}\subseteq G$ be a conjugacy class, and assume that for all $\sigma \in \mathcal{C}$, $\sigma\vert_{\mathbb{F}_{q^{n}}}=\Frob\vert_{\mathbb{F}_{q^{n}}}$.
    If
    \[q\geq 1336^{2}36^{d}N^{2\max\{d,3\}+6}m^{2\max\{d,3\}+6}n^{d+3}(\log(N))^{d+3}\ ,\]
    then $\pi_{Y/X}(1,\mathcal{C})>0$.
\end{Prop}
\begin{proof}
Denote by $U_{Y}\subseteq Y, U_{X}\subseteq X$ the open subsets of all unramified points in $Y$ and $X$ respectively, and take $\tilde{U}_{X},\tilde{U}_{Y}$ as in Proposition~\ref{Prop:4.1}.

From Theorem~\ref{Thm:1.3},
\[
    \bigg|\pi_{\tilde{U}_{Y}/\tilde{U}_{X}}(1,\mathcal{C})-\frac{\lvert \mathcal{C}\rvert}{m}q^{d}\bigg|\leq \dfrac{\lvert\mathcal{C}\rvert}{m}\left[(\delta-1)(\delta-2)q^{d-\frac{1}{2}}+24(2D+3)^{d+3}q^{d-1}\right]\, .
\]
Notice that if 
\[
    \dfrac{\lvert \mathcal{C}\rvert}{m}q^{d}>
    \dfrac{\lvert\mathcal{C}\rvert}{m}\left[(\delta-1)(\delta-2)q^{d-\frac{1}{2}}+24(2D+3)^{d+3}q^{d-1}\right]\, ,
\]
then $\pi_{Y/X}(1,\mathcal{C})\geq \pi_{U_{Y}/U_{X}}(1,\mathcal{C})>0$.
It is enough to prove that:
\[
    q>(\delta-1)(\delta-2)q^{1/2}+24(2D+3)^{d+3}\, ,
\]
meaning that we need to show that
\[q^{1/2}>\dfrac{(\delta-1)(\delta-2)+\sqrt{(\delta-1)^{2}(\delta-2)^{2}+96(2D+3)^{d+3}}}{2}\, ,\]
which can be checked to follow from
\[q\geq1336^{2}36^{d}N^{2\max\{d,3\}+6}m^{2\max\{d,3\}+6}n^{d+3}(\log(N))^{d+3}\, .\]
\end{proof}
\subsection{Proof of Theorem~\ref{thm:1App}}
We continue using Setting~\ref{Setting:app}.
\begin{Prop}\label{Prop:Gmax}
    Assume that $\Gal(L\overline{\mathbb{F}_{q}}/K\overline{\mathbb{F}_{q}})\cong S_{\deg_{t}(f)}$ and \[q>1336^{2}36^{d}N^{2\max\{d,3\}+6}m^{2\max\{d,3\}+6}n^{d+3}(\log(N))^{d+3}\, .\]
    Then, for every factorisation type $\lambda$ of a polynomial of degree $\deg_{t}(f)$, there exists a $\alpha\in \mathbb{F}_{q}^{u}$, such that $f(\alpha,t)$ is of factorisation type $\lambda$.
\end{Prop}
\begin{proof}
    Let $\lambda$ be a factorisation type of a polynomial of degree $\deg_{t}(f)$, and denote by $\mathcal{C}_{\lambda}\subseteq S_{\deg_{t}(f)}$ the conjugacy class of permutations with cycle type $\lambda$. We will first show that for $\alpha\in \mathbb{F}_{q}^{u}$, if 
    \[x_{\alpha}:=(x_{i}-\alpha_{i}\mid i\in [d])\in \tilde{U}_{X}\, ,\]
    where $\tilde{U}_{X}$ is as in Proposition~\ref{Prop:4.1}, then $f(\alpha,t)$ is of factorisation type $\lambda$ if and only if the Frobenius class of $x_{\alpha}$ is $\mathcal{C}_{\lambda}$.

    Indeed, let $y\in h^{-1}(x_{\alpha})$. Then, $\kappa(y)$ is the splitting field of $f(\alpha,t)$ over $\kappa(x_{\alpha})$. Therefore, the cycle type of
    \[\Frob_{\N x_{\alpha}}\in \Gal(\kappa(y)/\kappa(x_{\alpha}))\]
    as a permutations on the roots of $f(\alpha,t)$ is the factorisation type of $f(\alpha,t)$. Due to Corollary~\ref{Cor:e=1},
    \[D_{y/x}\cong \Gal(\kappa(y)/\kappa(x_{\alpha}))\, ,\]
    and so $\Frob_{\N x_{\alpha}}$ and $\Fr_{Y/X}[y/x_{\alpha}]$ have the same cycle type. Therefore, we get the desired equivalence.

    Thus, it is enough to show that $\pi_{\mathcal{C}_{\lambda}}(1)>0$. Since $q>2214^{2}6^{2d}N^{2\max\{d,3\}+6}m^{2\max\{d,3\}+6}(\log(N))^{d+3}$, this follows immediately from Proposition~\ref{Prop:PolyCheb}.
\end{proof}
Theorem~\ref{thm:1App} follows from this Proposition:
\begin{proof}
    In a similar fashion to \cite[Theorem 4.4.1]{Serre}, $\Gal(L\overline{\mathbb{F}_{q}}/K\overline{\mathbb{F}_{q}})\cong S_{\deg_{t}(f)}$, and the theorem follows from Proposition~\ref{Prop:Gmax}.
\end{proof}
\subsection{Explicit Bounds on The Bateman-Horn Conjecture Over Large Finite Fields}
Before we proceed with this application, we first prove a corollary of Theorem~\ref{Thm:1.2} that will prove to be useful.
\begin{Lemma}\label{Cor:app}
Let $\mathcal{C}\subseteq G$ be a conjugacy class, and assume that for all $\sigma\in \mathcal{C}$, $\sigma|_{\mathbb{F}_{q^{n}}}=\Frob_{q}|_{\mathbb{F}_{q^{n}}}$. Then,
\begin{align*}
    \left| \lvert\{\mathbf{a}\in \mathbb{F}_{q}^{d}\mid (x_{i}-a_{i}\mid i\in [d]) \text{ is unramified in } L\text{, }\Fr_{L/K}(\mathbf{a})=\mathcal{C}\}\rvert-\dfrac{\lvert\mathcal{C}\rvert}{m}q^{d}\right|&\\
    \leq
    \dfrac{\lvert\mathcal{C}\rvert}{m}\left[(\delta-1)(\delta-2)q^{d-\frac{1}{2}}+24(2D+3)^{d+3}q^{d-1}\right]+(D-1)q^{d-1}&\, .
\end{align*}
\end{Lemma}
\begin{proof}
    This is immediate from Theorem~\ref{Thm:1.3} and Proposition~\ref{Prop:4.1}.
\end{proof}
\begin{Not}
    Let $f\in \mathbb{F}_{q}[t]$ be a separable polynomial. Decompose $f$ into prime polynomials:
    \[f=\prod g_{i}\, .\]
    Denote by $\Theta_{f}\subseteq S_{\deg(f)}$ the conjugacy class of permutations with cycle type $(\deg(g_{i}))$. For multiple separable polynomials $f_{1},\ldots,f_{l}$, denote $\Theta_{f_{1},\ldots,f_{l}}=\Theta_{f_{1}}\times\ldots\times \Theta_{f_{l}}$.
\end{Not}
\begin{Corollary}\label{Cor:7}
Let $\mathcal{F}_{1},\ldots, \mathcal{F}_{l}\in \mathbb{F}_{q}[t][x],\ \deg_{x}(\mathcal{F}_{i})>0$ be distinct, monic in $x$, irreducible, and separable polynomials over $\mathbb{F}_{q}[t]$.
Let $d\geq \max_{i\in [l]}\{2,\Sl\mathcal{F}_{i}\}$ be a natural number, $a_{0},\ldots a_{d-1}$ be free variables,\ $\mathbf{a}=(a_{0},\ldots, a_{d-1})$ and \[f=t^{d}+\sum_{j<d} a_{j}t^{j},\ N_{i}=\deg_{t} \mathcal{F}_{i}(t,f)\, .\]
Let $L'$ be the splitting field of $\prod_{i\in [l]}\mathcal{F}_{i}(t,f)$ over $K:=\mathbb{F}_{q}(\mathbf{a})$. Then, $L'/K$ is a Galois extension. Let $G=\Gal(L'/K)$ be the corresponding Galois group. Let $\mathbb{F}_{q^{\nu}}=L'\cap \overline{\mathbb{F}_{q}}$, and denote
\[G_{1}=\{\sigma\in G\mid \sigma|_{\mathbb{F}_{q^{\nu}}}=\Frob_{q}\}\, .\]
We will think of $G$ as a subgroup of $S_{N_{1}}\times\ldots\times S_{N_{l}}$.

Then, for every conjugacy class $\mathcal{C}\subseteq S_{N_{1}}\times \ldots \times S_{N_{l}}$, we have:
\begin{align*}
    \Big|\lvert \{\tilde{f}\in \mathbb{F}_{q}[t]_{d}\mid \tilde{f} \text{ is monic},\ \Theta_{\mathcal{F}_{1}(t,\tilde{f}),\ldots,\mathcal{F}_{l}(t,\tilde{f})}=\mathcal{C}\}\rvert&-\dfrac{\lvert \mathcal{C}\cap G_{1}\rvert}{\lvert G_{1}\rvert}q^{d}\Big|\\
    &\leq\dfrac{\lvert\mathcal{C}\cap G_{1}\rvert}{\lvert G_{1}\rvert}q^{d}\left[(\delta-1)(\delta-2)q^{-\frac{1}{2}}+24(2D+3)^{d+3}q^{-1}\right]\\
    &\hphantom{.}+ (D-1+4N^{2}d)q^{d-1}\, ,
\end{align*}
where $D=(8N\log(N)+\nu)N\lvert G_{1}\rvert(2\lvert G_{1}\rvert-1)-3$, $\delta=DN\lvert G_{1}\rvert$, and $N=\deg(\prod_{i=1}^{l}\mathcal{F}_{i}(t,f))$.
\end{Corollary}
\begin{proof}
    Denote $\tilde{\mathcal{F}}_{i}(t,x)=\mathcal{F}_{i}(t,x+t^{d})$. Then, $L'$ is the splitting field of $\prod_{i\in [l]} \tilde{\mathcal{F}}_{i}(t,f-t^{d})$, and from~\cite[Lemma 2.1]{Entin}, $L'/K$ is Galois. For each $i\in [l]$, take
    \[\mathcal{G}_{i}=\mathcal{F}_{i}(t,f)\in \mathbb{F}_{q}[a_{0},\ldots,a_{d-1}][t]\]
    and
    \[g=\prod_{i\in [l]} \mathcal{G}_{i}\, .\]
    Then, $L=L'$, $\nu=n$. If we show that for any unramified ideal $(a_{i}-b_{i}\mid i<d)\lhd \mathbb{F}_{q}[a_{0},\ldots,a_{d-1}]$ for which all $\mathcal{F}_{i}(t,f(\textbf{b}))$ are separable,
    \[\Fr_{L/K}(\mathbf{b})\subseteq\Theta_{\mathcal{F}_{1}(t,f(\textbf{b})),\ldots,\mathcal{F}_{l}(t,f(\textbf{b}))}\cap G_{1}\, ,\]
    then from Lemma~\ref{Cor:app}, Proposition~\ref{Prop:4.1} and \cite[Theorem 1.2]{Rudnick} we will acquire the desired inequality.

    Indeed, if $\tilde{N}=\deg_{t}(g)$, $\alpha_{1},\ldots,\alpha_{\tilde{N}}$ are the roots of $g$ over $K$, $B=A[\alpha_{1},\ldots,\alpha_{\tilde{N}}]$ and $Q\lhd B$ is a maximal ideal lying over $(a_{i}-b_{i}\mid i<d)$, then
    \[\alpha_{i}+Q\in B/Q\]
    are the roots of
    \[g(b_{0},\ldots,b_{d-1},t)\, .\]
    Therefore, if
    \[\mathcal{G}_{i}(b_{0},\ldots,b_{d-1},t)=\prod h_{j,i}\]
    is the decomposition of the specialisation of $\mathcal{G}_{i}$ for some $i\in [l]$, then
    \[\Fr_{L/K}[Q/(a_{i}-b_{i}\mid i<d)] \mod Q\]
    acts as a cycle on the roots of $h_{j,i}$. All $h_{j,i}$'s are different, meaning that the cycles are disjoint. Therefore, for a cycle type $\lambda$, $\Fr_{L/K}(\underline{\mathbf{b}})=\mathcal{C}_{\lambda}$ if and only if $\Theta_{\mathcal{G}_{1}(b_{0},\ldots,b_{d-1},t),\ldots,\mathcal{G}_{l}(b_{0},\ldots,b_{d-1},t)}=\lambda$.
\end{proof}
We are now ready to prove Theorem~\ref{BatemanHorn}:
\begin{proof}
    Since \cite[Proposition 6.1]{Entin2} applies to
    \[\tilde{\mathcal{F}}_{i}(t,x):=\mathcal{F}_{i}(t,x+t^{d})\, ,\]
    we conclude that \cite[Proposition 2.4, Proposition 2.3]{Entin} hold true for $f=t^{d}+\sum_{j<d} a_{j}t^{j}$. Thus, the proof of \cite[Theorem 1.1]{Entin} shows that the error term of Corollary~\ref{Cor:7} is the same as the error term of $\lvert \pi_{q,\mathcal{F}_1,\ldots, \mathcal{F}_l}(d)-\rho^{-1}q^{d}\rvert$. Using the inequality $4N^{2}d\leq D-1$, we get the error term found in Theorem~\ref{BatemanHorn}.
\end{proof}
\begin{Corollary}
With the same notation as in Theorem~\ref{BatemanHorn}, if 
\[q>1336^{2}36^{d}N^{N(2\max\{d,3\}+6)}(\log(N))^{d+3}\, ,\] then
\[\{f\in \mathbb{F}_{q}[t]\mid \deg f=d, f\text{ is monic}.\ \forall{i\in [l]}, \mathcal{F}_{i}(t,f) \text{ is irreducible}\}\neq \emptyset\]
\end{Corollary}
\begin{proof}
By reducing to open subsets such as in Proposition~\ref{Prop:4.1}, this follows from Theorem~\ref{BatemanHorn} and the calculations done in Proposition~\ref{Prop:PolyCheb}.
\end{proof}
\section{Appendix}\label{Appendix}
We present a proof of Proposition~\ref{Prop:4.1}.

Let $f=\prod_{i\in [l]} f_{i}$ be the decomposition of $f$ into irreducible polynomials $f_{i}\in \mathbb{F}_{q}[x_{1},\ldots,x_{d}][t]$. Denote $N_{i}=\deg_{t}(f_{i})$. Without loss of generality, $N_{l}=\max\{N_{i}\mid i\in [l]\}$.
\begin{Not}
    For the next two lemmata, we will make use of the following notations:
For each $i\in [l], j\in [N_{i}]$, denote \[[(i,j)]=\{(i',j')\mid i'\in [l], j'\in [N_{i'}],(i',j')\leq_{\lex} (i,j)\}\, ,\]
where $\leq_{\lex}$ is the lexicographic order.

For every $(i,j)$, denote by $(i,j)+1$ the successor of $(i,j)$ in the lexicographic order.
\end{Not}
For every $i\in [l]$, let $\alpha_{i,1},\ldots,\alpha_{i,N_{i}}\in L$ be the roots of $f_{i}$ over $K$.
The following lemma shows that we can find primitive elements $\beta$ for $L\overline{\mathbb{F}_{q}}/K\overline{\mathbb{F}_{q}}$ of a "nice" form:
\[\beta=\sum c_{i,j}\alpha_{i,j},\ c_{i,j}\in \overline{\mathbb{F}_{q}} \]
such that any permutation on the $c_{i,j}$ will still yield a primitive element, and we can guarantee that most of these permutations will give us different elements.
\begin{remark}
    Notice that $S_{N_{1}}\times\ldots\times S_{N_{l}}$ acts on $[(l,N_{l})]$, via $(\tau_{1},\ldots,\tau_{l})(i,j)=(i,\tau_{i}(j))$, so we can embed $S_{N_{1}}\times\ldots\times S_{N_{l}}\subseteq S_{[(l,N_{l})]}$. 
\end{remark}
\begin{Lemma}\label{Lemma:5.1}
    Let $s=n\lceil\frac{4-2N+3N\log(N)+6\log(N)}{n}\rceil$. For each $i\in [l]$ there are $c_{i,1},\ldots,c_{i,N_{i}}\in \mathbb{F}_{q^{s}}$ such that for all $\tau\in S_{[(l,N_{l})]}$,
    \[\beta_{\tau}:=\sum_{(i,j)\in [(l,N_{l})]}c_{\tau(i,j)}\alpha_{i,j}\]
    is a primitive element for $L\mathbb{F}_{q^{s}}/K\mathbb{F}_{q^{s}}$, and for every two permutations $\sigma,\tau\in S_{[(l,N_{l})]}$,
    \[\beta_{\sigma}=\beta_{\tau}\iff \forall{(i,j)\in [(l,N_{l})]}, \alpha_{\sigma(i,j)}=\alpha_{\tau(i,j)}\] 
\end{Lemma}
\begin{proof}
    We shall prove a stronger statement than the one stated in the Lemma. We will prove that there are
    \[(c_{i,j})_{(i,j)\in [(l,N_{l})]} \in \mathbb{F}_{q^{s}}\]
    that satisfy the following:
    \begin{enumerate}[label=\arabic*.]
        \item For every $\emptyset \neq \Lambda\subseteq [(l,N_{l})]$ and every injection $\tau\colon \Lambda\to [(l,N_{l})]$, the element \[\beta_{\Lambda,\tau}:=\sum_{\omega\in \Lambda} c_{\omega}\alpha_{\tau(\omega)}\]
    is a primitive element of $K\overline{\mathbb{F}_{q}}(\alpha_{\tau(\omega)}\mid \omega\in \Lambda)/K\overline{\mathbb{F}_{q}}$.
    \item For every $\emptyset\neq \Lambda\subseteq [(l,N_{l})]$ and $\sigma,\tau\colon \Lambda\to [(l,N_{l})]$ injections, if for some $\omega\in \Lambda$, $\alpha_{\sigma(\omega)}\neq \alpha_{\tau(\omega)}$, then $\beta_{\Lambda,\sigma}\neq \beta_{\Lambda,\tau}$.
    \end{enumerate}
    We will choose $c_{i,j}$ inductively, with respect to the lexicographic order on $[(l,N_{l})]$.
    
    Take $c_{1,1}=1$.
    
    Let $i_{0}\in [l], j_{0}\in [N_{i_{0}}]$, and assume that we chose
    \[c_{1,1},\ldots,c_{i_{0},j_{0}}\]
    in such a way that $1$ and $2$ are satisfied whenever we limit ourselves to $\Lambda\subseteq [(i_{0},j_{0})]$.
    
    From the proof of the primitive element theorem, we know that if $E$ is an infinite field, and $E(\gamma,\delta)/E$ is a finite separable extension, then every
    \[c\in E\setminus \{\dfrac{\sigma(\gamma)-\gamma}{\delta-\sigma(\delta)}\mid \sigma\in \Aut(E(\gamma,\delta)/E) \text{ does not fix }\delta\}\]
    satisfies that $\gamma+c\delta$ is a primitive element.
    We will use this criterion to guarantee the primitiveness of our elements.
    
    Informally, this condition will guarantee that up to a finite number of elements, every choice of $c_{(i_{0},j_{0})+1}$ will suffice. The formalism is riddled with indices, but we provide it nonetheless:

    For every subset $\emptyset\neq \Lambda\subseteq [(i_{0},j_{0})]$, every injection $\tau\colon \Lambda\to [(l,N_{l})]$, every element
    \[\alpha\in \{\alpha_{i,j}\mid (i,j)\in [(l,N_{l})]\setminus \tau[I]\}\, ,\]
    and every automorphism $\sigma\in \Gal(L\mathbb{F}_{q^{s}}/K\mathbb{F}_{q^{s}})$ that does not fix $\alpha$, denote:
    \[A_{\Lambda,\tau,\alpha,\sigma}=\dfrac{\sum_{\omega\in \Lambda} c_{\omega}(\sigma(\alpha_{\tau(\omega)})-\alpha_{\tau(\omega)})}{\alpha-\sigma(\alpha)}\]
    Additionally, for every
    $\Lambda\subseteq [(i_{0},j_{0})]$, every injections $\sigma,\tau\colon \Lambda\to [(l,N_{l})]$ for which there exists an $\omega\in \Lambda$ such that $\alpha_{\sigma(\omega)}\neq \alpha_{\tau(\omega)}$, and every two elements \[\alpha\in \{\alpha_{i,j}\mid (i,j)\in [(l,N_{l})]\setminus \sigma[\Lambda]\},\ \alpha\neq \alpha'\in \{\alpha_{i,j}\mid (i,j)\in [(l,N_{l})]\setminus \tau[\Lambda]\}\]
    denote
    \[B_{\Lambda,\sigma,\tau,\alpha,\alpha'}=\dfrac{\sum_{\omega\in \Lambda}c_{\omega}\alpha_{\tau(\omega)}-\sum_{\omega\in \Lambda}c_{\omega}\alpha_{\sigma(\omega)}}{\alpha-\alpha'}\]
    And take 
    \begin{multline*}
        \mathcal{A}=\{A_{\Lambda,\tau,\alpha,\sigma}\mid \Lambda\subseteq [(i_{0},j_{0})], \tau\colon \Lambda\to [(l,N_{l})]\text{ an injection},
        \alpha\in \{\alpha_{(i,j)}\mid (i,j)\in [(l,N_{l})]\setminus \tau[\Lambda]\},\\ \sigma\in \Gal(L\mathbb{F}_{q^{s}}/K\mathbb{F}_{q^{s}}) \text{ does not fix }\alpha\}\, ,
    \end{multline*}
    and
    \begin{multline*}
        \mathcal{B}=\{B_{\Lambda,\sigma,\tau,\alpha,\alpha'}\mid \Lambda\subseteq [(i_{0},j_{0})],\sigma,\tau\colon \Lambda\to [(l,N_{l})] \text{ injections, }
        \alpha\in \{\alpha_{(i,j)}\mid (i,j)\in [(l,N_{l})]\setminus \sigma[\Lambda]\},\\ \alpha'\in \{\alpha_{(i,j)}\mid (i,j)\in [(l,N_{l})]\setminus \tau[\Lambda]\}\setminus \{\alpha\}\}\, .
    \end{multline*}
    Since
    \[|\mathcal{A}\cup \mathcal{B}|\leq 2(2^{N}-1)N!^{2}(N-1)^{2}\, ,\]
    we can take $c_{(i,j)+1}\in \mathbb{F}_{q^{s}}\setminus (\mathcal{A}\cup \mathcal{B})$.

    $\mathcal{B}$ guarantees that the function is injective. From the condition that implies primitivity, we achieve the primitivity of all the desired elements.
\end{proof}
\begin{remark}
    Notice that we can embed
    \[\Gal(L\overline{\mathbb{F}_{q}}/K\overline{\mathbb{F}_{q}})\hookrightarrow S_{N_{1}}\times\ldots\times S_{N_{l}}\subseteq S_{[(l,N_{l})]}\, .\]
    Take $\sigma\in \Gal(L\overline{\mathbb{F}_{q}}/K\overline{\mathbb{F}_{q}})$, and let it act on $[(l,N_{l})]$ in the following way:
    
    $\sigma(i,j)\in \{i\}\times [N_{i}]$ is the unique element in that set satisfying $\sigma(\alpha_{i,j})=\alpha_{\sigma(i,j)}$.
\end{remark}
The next lemma will show that we can find a primitive element for $L\overline{\mathbb{F}_{q}}/K\overline{\mathbb{F}_{q}}$ for which we can bound the total degree of its minimal polynomial.
\begin{Lemma}\label{Lemma:5.2}
    Let $s=n\lceil\frac{4-2N+3N\log(N)+6\log(N)}{n}\rceil$. There exist $(c_{i,j})_{(i,j)\in [(l,N_{l})]}\in \mathbb{F}_{q^{s}}$, such that
    \[\beta:=\sum_{(i,j)\in [(l,N_{l})]}c_{i,j}\alpha_{i,j}\]
    is a primitive element for $L\mathbb{F}_{q^{s}}/K\mathbb{F}_{q^{s}}$, and its minimal polynomial $h\in \mathbb{F}_{q^{s}}[x_{1},\ldots,x_{d}][t]$ satisfies
    \[\deg(h)\leq Nm\, .\]
\end{Lemma}
\begin{proof}
    Take $(c_{i,j})_{(i,j)\in [(l,N_{l})]}\in \mathbb{F}_{q^{s}}$ as described in Lemma~\ref{Lemma:5.1}. For every permutation ${\tau\in S_{[(l,N_{l})]}}$, denote
    \[\beta_{\tau}=\sum_{(i,j)\in [(l,N_{l})]}c_{i,j}\alpha_{\tau(i,j)}\]
    and let $h_{\tau}\in \mathbb{F}_{q^{s}}[x_{1},\ldots,x_{d}][t]$ be its minimal polynomial over $K\mathbb{F}_{q^{s}}$. We will show that there exists a $\tau\in S_{[(l,N_{l})]}$ such that
    \[\deg(h_{\tau})\leq Nm\, .\]
    Define an equivalence relation on $S_{[(l,N_{l})]}$:
    \[\tau\equiv \tilde{\tau}\iff \exists\sigma\in \Gal(L\overline{\mathbb{F}_{q}}/K\overline{\mathbb{F}_{q}}),\ \sigma(\beta_{\tau})=\beta_{\tilde{\tau}}\]
    Denote
    \[u=\{\sigma\in S_{[(l,N_{l})]}\mid \forall{(i,j)\in [(l,N_{l})]}, \alpha_{\sigma(i,j)}=\alpha_{(i,j)}\},\ \tilde{N}=\sum_{i\in [l]}N_{i}=\deg_{t}(f)\, .\]
    We will first show that there are exactly $\frac{\tilde{N}!}{um}$ equivalence classes. Let $\tau\in S_{[(l,N_{l})]}$ be a permutation.
    Let
    \[\{\sigma_{1},\ldots,\sigma_{m}\}\subseteq S_{N_{1}}\times\ldots\times S_{N_{l}}\subseteq S_{[(l,N_{l})]}\]
    be the image of $\Gal(L\overline{\mathbb{F}_{q}}/K\overline{\mathbb{F}_{q}})$ in $S_{N_{1}}\times\ldots\times S_{N_{l}}$.
    
    Then, since for every $i\in [m]$, \[\beta_{\sigma_{i}\circ\tau}=\sigma_{i}(\beta_{\tau})\, ,\]
    we see that
    \[\{\beta_{\tilde{\tau}}\mid \tau\equiv \tilde{\tau}\}=\{\beta_{\sigma_{i}\circ\tau}\mid i\in [m]\}\, .\]
    Additionally, using Lemma~\ref{Lemma:5.1}, we see that \[u=\{\tilde{\tau}\in S_{[(l,N_{l})]}\mid \beta_{\tau}=\beta_{\tilde{\tau}}\}\, ,\]
    and that for every $i\neq j$,
    \[\beta_{\sigma_{i}\circ \tau}\neq \beta_{\sigma_{j}\circ \tau}\, .\]
    Therefore, there are exactly $um$ elements in the equivalence class of $\tau$. In particular, there are $\frac{\tilde{N}!}{um}$ equivalence classes.
    
    Let ${\tau}_{1},\ldots,{\tau}_{\frac{\tilde{N}!}{um}}\in S_{[(l,N_{l})]}$ be representatives for this equivalence relation.
    Notice that
    \[\prod_{v=1}^{\tilde{N}!/um} h_{\tau_{v}}^{s}=\prod_{\tau\in S_{[(l,N_{l})]}}(t-\sum_{(i,j)\in [(l,N_{l})]} c_{i,j}\alpha_{\tau(i,j)})\, .\]
    To achieve a bound on $\deg(h_{\tau_{v}})$ for some $v$, we will first bound the degree of 
    \[\tilde{h}:=\prod_{\tau\in S_{[(l,N_{l})]}}(t-\sum_{(i,j)\in [(l,N_{l})]} c_{i,j}\alpha_{\tau(i,j)})\, .\]
    Let $v\leq \tilde{N}!$, and define:
    \[P_{v}(t_{1,1},\ldots,t_{l,N_{l}})=\Sym^{\tilde{N}!}_{v}(\sum_{(i,j)\in [(l,N_{l})]}c_{i,j}t_{\tau(i,j)}\mid \tau\in S_{[(l,N_{l})]})\, ,\]
    where $\Sym^{\tilde{N}!}_{v}$ is the elementary symmetric polynomial of $\tilde{N}!$ variables, of degree $v$.
    This is a symmetric polynomial, so there exists a polynomial $Q_{v}\in \overline{\mathbb{F}_{q}}[y_{1},\ldots,y_{N}]$ such that
    \[P_{v}=Q_{v}\circ(\Sym^{\tilde{N}}_{1},\ldots,\Sym^{\tilde{N}}_{\tilde{N}})\, .\]
    Since $P_{v}$ is homogeneous of degree $v$, we can take $Q_{v}$ in such a way that for any monomial $M$ appearing in $Q_{v}$,
    \[M\circ (\Sym^{\tilde{N}}_{1},\ldots,\Sym^{\tilde{N}}_{\tilde{N}})\]
    is homogeneous of degree $v$.
    Therefore, since $\Sym^{\tilde{N}}_{i}$ is the coefficient of $t^{\tilde{N}-i}$ in $f$: 
    \[\deg(P_{v}(\alpha_{1,1},\ldots\alpha_{l,N_{l}}))\leq v\max_{i\in [\tilde{N}]}\deg(\Sym^{\tilde{N}}_{i}(\alpha_{1,1},\ldots,\alpha_{l,N_{l}}))\leq Nv\, .\]
    $P_{v}(\alpha_{1,1},\ldots,\alpha_{l,N_{l}})$ is the coefficient of $t^{\tilde{N}!-v}$ in $\tilde{h}$, so we conclude that
    \[\deg(\tilde{h})\leq \max_{v\leq \tilde{N}!}\{Nv+\tilde{N}!-v\}=\max_{v\leq \tilde{N}!}\{(N-1)v+\tilde{N}!\}= (N-1)\tilde{N}!+\tilde{N}!=N\cdot \tilde{N}!\, .\]
    Therefore, if we order $h_{\tau_{1}},\ldots h_{\tau_{\frac{N!}{um}}}$ in such a way that $\deg(h_{\tau_{v}})\leq \deg(h_{\tau_{v+1}})$, we can see that:
    \[\dfrac{\tilde{N}!}{um}\deg(h_{\tau_{1}}^{u})\leq \deg(\tilde{h})\leq N\cdot \tilde{N}!\, .\]
    Therefore:
    \[\deg(h_{\tau_{1}})\leq Nm\, .\]
\end{proof}
We can now prove Proposition~\ref{Prop:4.1}:
\begin{Prop*}
    Let $U_{Y}\subseteq Y, U_{X}\subseteq X$ be the open subsets of all unramified points in $Y$ and $X$ respectively. Then, there exist open subsets $\tilde{U}_{X}\subseteq U_{X}$, $\tilde{U}_{Y}\subseteq U_{Y}$ such that $\tilde{U}_{Y}=h^{-1}[\tilde{U}_{X}]$, and: 
    \begin{enumerate}[label=\arabic*.]
        \item There exists an embedding of ${\tilde{U}_{Y}{}}_{\overline{\mathbb{F}_{q}}}$ into the affine space $\mathbb{A}^{d+2}_{\overline{\mathbb{F}_{q}}}$, such that ${\tilde{U}_{Y}{}}_{\overline{\mathbb{F}_{q}}}$ is defined by the vanishing of $2$ polynomials of degrees bounded by $D$, and $\delta_{{\tilde{U}_{Y}{}}_{\overline{\mathbb{F}_{q}}}}\leq \delta$.
        \item \[\lvert U_{X}\setminus \tilde{U}_{X}\rvert_{1}\leq (D-1)q^{d-1}\, .\]
    \end{enumerate}
\end{Prop*}
\begin{proof}
    Let $\beta\in L\mathbb{F}_{q^{s}}$ be a primitive element for $L\mathbb{F}_{q^{s}}/K\mathbb{F}_{q^{s}}$ as described in Lemma~\ref{Lemma:5.2}, and denote by $h\in \overline{\mathbb{F}_{q}}[x_{1},\ldots,x_{d}][t]$ its minimal polynomial over $K\overline{\mathbb{F}_{q}}$. 
    
    Take 
    \[\tilde{U}_{X}=X\cap \{\Norm_{K\mathbb{F}_{q^{s}}/K}(\Disc_{t}(h))\neq 0\} =U_{X}\cap \{\Norm_{K\mathbb{F}_{q^{s}}/K}(\Disc_{t}(h))\neq 0\}\, ,\]
    and
    \[\tilde{U}_{Y}=Y\cap\{\Norm_{K\mathbb{F}_{q^{s}}/K}(\Disc_{t}(h))\neq 0\} =U_{Y}\cap \{\Norm_{K\mathbb{F}_{q^{s}}/K}(\Disc_{t}(h))\neq 0\}=h^{-1}[\tilde{U}_{X}]\, .\]
    Then, for $(1)$,
    \begin{align*}
    (\tilde{U}_{Y})_{\overline{\mathbb{F}_{q}}}&\cong \Spec \overline{\mathbb{F}_{q}}[x_{1},\ldots,x_{d}][y][\beta]/(y\Norm_{K\mathbb{F}_{q^{s}}/K}(\Disc_{t}(h))-1)\\
    &=\Spec \overline{\mathbb{F}_{q}}[x_{1},\ldots,x_{d}][t,y]/(h,y\Norm_{K\mathbb{F}_{q^{s}}/K}(\Disc_{t}(h))-1)\, .
    \end{align*}
    so we can embed $(\tilde{U}_{Y})_{\overline{\mathbb{F}_{q}}}$ into a $d+2$-dimensional affine space such that the image is defined by the vanishing of $2$ polynomials, of degrees bounded by:
    \begin{align*}
    \max\{\deg(h), \deg(y\Norm_{K\mathbb{F}_{q^{s}}/K}(\Disc_{t}(h))-1)\}&\leq\max\{Nm, (8N\log(N)+n-N)Nm(2m-1)+1\}\\
    &\leq (8N\log(N)+n)Nm(2m-1)-3=D\, ,
    \end{align*}
    where we used the inequlaity $s\leq 8N\log(N)+n-N$. Additionally,
    \[\delta_{{\tilde{U}_{Y}{}}_{\overline{\mathbb{F}_{q}}}}\leq \deg(h)\deg(y\Norm_{K\mathbb{F}_{q^{s}}/K}(\Disc_{t}(h))-1)\leq Nm((8N\log(N)+n)Nm(2m-1)-3)\, =\delta.\]

    For $(2)$,
    $U_{X}\setminus \tilde{U}_{X}=U_{X}\cap \{\Norm_{K\mathbb{F}_{q^{s}}/K}(\Disc_{t}(h))=0\}$, and therefore,
    \begin{align*}
    \lvert U_{X}\setminus \tilde{U}_{X}\rvert_{1}&\leq \lvert \Spec \mathbb{F}_{q}[x_{1},\ldots,x_{d}]/(\Norm_{K\mathbb{F}_{q^{s}}/K}(\Disc_{t}(h)))\rvert_{1}\\
    &\leq ((8N\log(N)+n)Nm(2m-1)-4)q^{d-1}\, .
    \end{align*}
\end{proof}
\printbibliography
\end{document}